\documentclass{amsart}

\usepackage[letterpaper, left=1.1in, right=1.1in, bottom=1in, top=1in]{geometry}

\usepackage{tikzsymbols}
\usepackage{textcomp}
\usepackage[dvipsnames]{xcolor}
\usepackage{fancyhdr}
\usepackage{graphicx}
\usepackage{empheq}
\usepackage{ifthen}
\usepackage{amsmath,amssymb,amsthm,mathrsfs}
\usepackage{gastex}
\usepackage{color,graphicx}
\usepackage{tikz-cd}
\usetikzlibrary{decorations.pathmorphing,shapes,matrix}
\usepackage{mdframed}
\usepackage{caption}
\usepackage{lipsum}
\usepackage{mathabx}
\usepackage{url}
\usepackage{subcaption}
\usepackage{outlines} 
\usepackage{pgfplots}
\usepackage{ytableau}
\pgfplotsset{compat=1.18}
\usepackage{thmtools}
\usepackage{hyperref}
\hypersetup{
    colorlinks,
    linkcolor={red!50!black},
    citecolor={blue!50!black},
    urlcolor={blue!80!black}
}
\usepackage{cleveref}
\usepackage{float}

\DeclareMathOperator{\Aut}{Aut}

\DeclareMathOperator{\sym}{Sym}

\DeclareMathOperator{\inv}{inv}

\DeclareMathOperator{\code}{{code}}

\DeclareMathOperator{\rank}{{rank}}

\DeclareMathOperator{\W}{\mathbf{W}}

\newcommand{\N}{{\mathbb N}}
\newcommand{\Z}{{\mathbb Z}}
\newcommand{\C}{{\mathbb C}}

\newcommand{\BC}{\leq_{BC}}

\newcommand{\bc}{<_{BC}}

\DeclarePairedDelimiter{\floor}{\lfloor}{\rfloor}
\DeclarePairedDelimiter{\ab}{\lvert}{\rvert}

\newcommand{\DeclareAutoPairedDelimiter}[3]{%
\expandafter\DeclarePairedDelimiter\csname Auto\string#1\endcsname{#2}{#3}%
  \begingroup\edef\x{\endgroup
    \noexpand\DeclareRobustCommand{\noexpand#1}{%
      \expandafter\noexpand\csname Auto\string#1\endcsname*}}%
  \x}

\DeclareAutoPairedDelimiter{\set}{\{}{\}}

\def\sus{\subseteq}
\def\la{\lambda}
\def\es{\varnothing}

\def\ol{\overline}
\def\x{\times}

\def\ox{\otimes}
\def\N{\mathbb{N}}
\def\C{\mathbb{C}}
\def\Z{\mathbb{Z}}

\def\D{\mathbb{D}}
\renewcommand{\P}{\mathbb{P}}

\def\mc#1{\mathcal{#1}}
\def\mb#1{\mathbf{#1}}

\definecolor{darkred}{rgb}{0.7,0,0} 
\newcommand{\defncolor}{\color{darkred}}
\newcommand{\dfn}[1]{{\defncolor\emph{#1}}} 
\newcommand{\ul}[1]{\underline{#1}}

\theoremstyle{plain}
\newtheorem{theorem}{Theorem}[section]
\newtheorem{lemma}[theorem]{Lemma}
\newtheorem{prop}[theorem]{Proposition}
\newtheorem{cor}[theorem]{Corollary}
\newtheorem{openp}[theorem]{Open Problem}

\theoremstyle{definition}
\newtheorem{defn}[theorem]{Definition}

\newtheorem{example}[theorem]{Example}
\newtheorem{remark}[theorem]{Remark}

\numberwithin{theorem}{section}
\numberwithin{equation}{section}

\usepackage[colorinlistoftodos]{todonotes}

\usepackage{setspace}

\date{\today}
\begin{document}
\title{On Symmetric Chain Decompositions for Bruhat--Chevalley Posets} 
\author[A. Bingham]{Aram Bingham}\address[A.~Bingham]{Departamento de Matem\'aticas, Universidad de Chile, 3425 Las Palmeras, Ñuñoa, Santiago, CL}
\email{\textcolor{blue}{aram@matmor.unam.mx}}

\author[N.~F.~Díaz Morera]{N. Fernando Díaz Morera}
\address[N.~F.~Díaz Morera]{Department of Mathematics, Fitchburg State University,
MA, USA}
\email{\textcolor{blue}{ndiazmor@fitchburgstate.edu}}

\author[Ö.~Uğurlu]{Özlem Uğurlu}
\address[Ö.~Uğurlu]{Department of Mathematics and Statistics,  Saint Louis University, MO, USA}
\email{\textcolor{blue}{ozlem.ugurlu@slu.edu}}
\begin{abstract}
   We consider Stanley's problem on describing symmetric chain decompositions for the Bruhat--Chevalley posets of partial flag varieties and, more generally, Schubert varieties. These posets describe the containment relations of cell closures in the partial flag variety $G/P$ of a connected complex reductive algebraic group $G$. A famous instance of this problem asks whether the poset of partitions that fit inside a rectangle ordered by containment admits a symmetric chain decomposition (SCD), corresponding to the case of Grassmannian Schubert varieties in type $A$. Our first main result applies the ``complete BP decomposition" of Richmond--Slofstra to show that the SCD problem for any rationally smooth Schubert variety reduces to the (generalized) Grassmannian case, in any semisimple type. We further study the SCDs provided by codes for Bruhat--Chevalley posets in type $A$. These SCDs can be viewed as a consequence of a particular form of complete BP decomposition. In particular, we show that a poset of permutation-labeled partitions in a rectangle corresponding to the quotient $\mathcal{S}_n/\mathcal{S}_{n-k}$ admits an SCD. Further, posets of weak compositions inside a rectangle are realized as Bruhat--Chevalley intervals with associated codes and SCDs.
\end{abstract}
\maketitle

\section{Introduction}

A finite ranked poset $\mc{P}$ is said to admit a \emph{symmetric chain decomposition} if there exists a collection of disjoint saturated chains $C_1,\dots,C_N$ that are symmetric about the middle rank such that $\mc{P}=\bigsqcup_{i=1}^N C_i$. Such a poset is called a \emph{symmetric chain order} (\emph{SCO}). The study of symmetric chain decompositions (\emph{SCDs}) originates in the proof by de Bruijn--van Ebbenhorst Tengbergen--Kruyswijk~\cite{debruijn_etal} that the poset of divisors of an integer satisfies the ``Sperner property,'' generalizing a result of Sperner~\cite{sperner} on the Boolean lattice. From the existence of an SCD, it follows that a poset is also rank-symmetric, rank-unimodal, and ``strongly Sperner.'' Following~\cite{proctorSaksSturtevant}, a poset with this trio of properties in conjunction is sometimes called \emph{Peck}, after the fictitious mathematician G. W. Peck conjured by Kleitman and collaborators~\cite{peck}.

Let $G$ be a connected complex reductive algebraic group with parabolic subgroup $P$ containing a fixed Borel subgroup $B$. Then the \emph{Bruhat--Chevalley order} (often known as ``strong Bruhat" or just ``Bruhat" order in the literature) describes the containment relations between $B$-orbit closures of the (partial) flag variety $G/P$. The $B$-orbits have the structure of affine cells often known as \emph{Schubert cells} and their closures are \emph{Schubert varieties}. These varieties, denoted $X^P(w):=\ol{BwP/P}$, are parameterized by coset representatives $w$ of the quotient $\W/\W_P$, where $\W$ is the Weyl group of $G$ and $\W_P$ is the Weyl group of $P$. Their containment relations then induce a Bruhat--Chevalley order relation $\BC$ on a set $\W^P$ of coset representatives for $\W/\W_P$. That is, for $v,w\in \W^P$,
\[v\BC w \iff X^P(v)\sus X^P(w).\]
The resulting poset $(\W^P, \BC)$ inherits many properties from the Bruhat--Chevalley order on $G/B$, captured by $(\W,\BC)$. A lower interval $[e,w]^P$ of $(\W^P,\BC)$ then describes a cell decomposition of the Schubert variety $X^P(w)$. Consequently, Bruhat--Chevalley intervals play a crucial role in much of combinatorial algebraic geometry and representation theory.

In~\cite{stanley1980}, Stanley used the hard Lefschetz theorem to prove that for smooth $X^P(w)$ (a class which includes the full posets $(\W^P, \BC)$), $[e,w]^P$ is Peck. He then posed the following open problem.\footnote{The problem is actually stated for posets obtained from a smooth, irreducible, complex projective variety with a cellular decomposition, with cells ordered by closure containment, but Bruhat--Chevalley posets are the primary examples of interest.}
\begin{openp}\label{prob:gen}
    Let $[e,w]^P$ be the Bruhat--Chevalley poset of Schubert varieties contained in a (rationally) smooth $X^P(w)$ inside the homogeneous space $G/P$. Is $[e,w]^P$ a symmetric chain order?
\end{openp}
Note that the hypothesis of smoothness was included in~\cite{stanley1980} to guarantee Poincaré duality, and hence rank-symmetry of the Bruhat--Chevalley poset. It was later shown by Carrell--Peterson~\cite{carrellPeterson} that for Schubert varieties, the weaker property of rational smoothness is sufficient to guarantee rank-symmetry. However, rational smoothness implies smoothness in simply-laced type~\cite{carrellKuttler}.

Stanley also observed \cite[p. 182]{stanley1980} that the Bruhat--Chevalley order on Weyl groups of type $A$, $B$, and $D$ are all SCOs. Since then, most of the work on~\Cref{prob:gen} has focused on the special case of $(\W^P, \BC)$ where $\W\cong \mc{S}_n$ and $\W_P\cong \mc{S}_k\x \mc{S}_{n-k}$, corresponding to the situation where $G/P$ is the \emph{Grassmannian variety} of $k$-planes inside of an $n$-dimensional complex vector space. It is well known that the Bruhat--Chevalley order on the Schubert varieties of the Grassmannian are in bijection with partitions that fit inside of a rectangle. Thus, we  have the following specialization of~\Cref{prob:gen}.

\begin{openp}\label{prob:rect}
    Does the poset $L(k,n-k)$ of partitions inside a $k \x (n-k)$ rectangle ordered by containment of shapes admit a symmetric chain decomposition for any $1\leq k<n$?
\end{openp}

 The origins of this problem can be traced back to the proof of the unimodality of $q$-binomial coefficients, first established by Sylvester~\cite{sylvester}. However, construction of an SCD for these posets has only been achieved for $k\leq 6$~\cite{riess, wen24, wen26}.

 In this article, we make progress on~\Cref{prob:gen} by showing that it can be reduced to a generalized Grassmannian case in any semisimple Lie type. In particular, our first main result (\Cref{thm:main1}) implies the following.
\begin{theorem}\label{thm:intro}
     In type $A$, $[e,w]^P$ is an SCO for all (rationally) smooth $(w,P)$ if and only if $L(k,n-k)$ is an SCO for all $1\leq k<n$. That is, \Cref{prob:gen} and \Cref{prob:rect} are equivalent in type $A$.
 \end{theorem}
 This result is a fairly simple consequence of work by Richmond--Slofstra \cite{richmondSlofstra}, establishing what are called \emph{complete BP decompositions} in the summary article \cite{ohRichmond}. These decompositions correspond to iterated fiber bundle structures on rationally smooth Schubert varieties. One consequence of such a decomposition is a ``factorization" of the the interval $[e,w]^P$, given as an order-preserving bijection from a product of Grassmannian intervals.  The standard result that the product of SCOs is again an SCO then completes the proof. 
 
 Furthermore, we establish SCDs for certain Bruhat--Chevalley posets which are ``halfway" between $\mc{S}_n$ and $\mc{S}_n/\mc{S}_k \x \mc{S}_{n-k}$. \Cref{thm:WJ_SCO} describes an order-preserving bijection from a product of chains to the Bruhat--Chevalley order on $\W/\W_P=\mc{S}_n/\mc{S}_{n-k}$ for any $k$ (see also \Cref{rk:Sn-k}). As products of chains are SCOs, so is $(\W^P,\BC)$.  This can be viewed as a consequence of the existence of what we call a \emph{g-code} (\Cref{def:gcode}) for such posets. A related result was shown independently by Gaetz--Gao \cite[Theorem 5.3]{gaetzGao}, who further relate the generalized notion of ``Lehmer" codes\footnote{This term usually applies to a bijection only for permutations of $\mc{S}_n$. Coupled with that it was described by Charles-Ange Laisant 
 \cite{laisant} over 70 years before D.H. Lehmer \cite{lehmer}, we find the latter's name inappropriate to our general setting.} from \cite{bologniniSentinelli2025} to the existence of certain BP decompositions for elements of a given $\W$. Under certain conventions, the g-code for $\mc{S}_n/\mc{S}_{n-k}$ is the restriction of the usual code for $\mc{S}_n$ to the set of shortest length coset representatives. A g-code for an interval $[e,w]^P$ in type $A$ can also be viewed as a complete BP decomposition where each Grassmannian factor is a projective space $\C\P^{j}$ for some $j$.
 
 Elements of the poset $\mc{S}_n/\mc{S}_{n-k}$ are naturally in correspondence with pairs consisting of a partition in a rectangle (i.e. an element of $\mc{S}_n/(\mc{S}_k \x \mc{S}_{n-k})$) and a permutation in $\mc{S}_k$. We give explicit details in \Cref{prop:coset_pairs}. This bijection is also order- and length-preserving. A more general version of this statement is given in \Cref{prop:cox_subquotients_bij} and holds for any Coxeter group $\W$ and subgroup $\W_K$ such that $\W_K=\W_I \x \W_J$. The bijective poset map $\W^K \x \W_I \to \W^J$ is another particular instance of a BP decomposition, though it is not generally a poset isomorphism. 

 Finally, we show that $(\mc{S}_n/\mc{S}_{n-k}, \BC)$ contains an interval $[e,w]^P$ of rank $k(n-k)$ which contains $L(k,n-k)$ as an equal-rank subposet. More generally, in \Cref{thm:interval_bij} we show that the product of $k$ chains of any equal-length $s$ is in order-preserving bijection with a Bruhat--Chevalley interval. Thus, these intervals also admit g-codes and SCDs, further raising the question of the existence of g-codes and SCDs for parabolic quotients and intervals. We plan to return to this question with some results in types $B$ and $D$ in forthcoming work.

Our paper is structured as follows. In \Cref{sec:background} we recall necessary background on Coxeter groups, Bruhat--Chevalley order, BP decompositions, codes, and symmetric chain decompositions. In \Cref{sec:results} we state and prove our main results described above. Finally, in \Cref{sec:future} we comment on the relation of our results to previous work on SCOs, rational smoothness of Schubert varieties, and factorization of Poincaré polynomials, and we mention some future directions.

\section{Background and notation} 
\label{sec:background}

\subsection{Bruhat--Chevalley order and parabolic quotients}

We follow the notation and conventions of ~\cite{bjornerBrenti} and~\cite{Humphreys}. Let $\Z^{+}$ denote the set of positive integers and let $n\in \Z^+$.
We define a \dfn{Coxeter system} to be a pair $(\W, S)$ consisting of a group $\W$ and a set of generators $S \subset \W$ (also called \emph{simple reflections}), subject only to relations of the form
$$
\left(s s^{\prime}\right)^{m\left(s, s^{\prime}\right)}=1,
$$
where $m(s, s)=1$ and $m\left(s, s^{\prime}\right)=m\left(s^{\prime}, s\right) \geq 2$ for $s \neq s^{\prime}$ in $S$. We will only consider finite Coxeter groups in this article. We say that a finite $(\W, S)$ is \dfn{crystallographic} if $m(s,s')\in\set{2,3,4,6}$. The finite crystallographic Coxeter groups are exactly those which appear as the Weyl group of a connected reductive algebraic group. The \dfn{Coxeter graph} is constructed as the graph with vertices corresponding to generators and edges $\set{s, s'}$ for each pair with $m(s,s')\geq 3$ labeled by the value $m(s,s')$. A Coxeter system is said to be \dfn{irreducible} if its Coxeter graph is connected.

Let $T= \set{wsw^{-1}\;\mid\; w\in \W, s\in S}$ denote the set of \emph{reflections} of $\W$. Recall that every $w\in \W$ can be written as a product of simple reflections $w=s_{i_1}s_{i_2}\cdots s_{i_k}$. If $k$ is minimal among all such expressions for $w$ then $k$ is called the \dfn{length} of $w$, written $\ell(w) = k$, and $s_{i_1}\cdots s_{i_k}$ is called a \dfn{reduced word} for $w$. We now define two subsets of $T$ which will be used to define descents of $w$ as follows. The set of left associated reflections to $w$ is defined by: 
\[T_L(w) = \set{t\in T \;\mid\; \ell(tw) < \ell(w)}.\]
The set of right associated reflections to $w$ is defined by: 
\[T_R(w) = \set{t\in T \;\mid\; \ell(wt) < \ell(w)}.\]
For $w\in \W$, the set $D_L(w) = T_L(w) \cap S$ is called the \dfn{left descent set} and the set $D_R(w) = T_R(w) \cap S$ is called the \dfn{right descent set}. By  \cite[Corollary 1.4.6]{bjornerBrenti}, $s\in D_R(w)$ if and only if some reduced expression for $w$ ends with the letter $s$. 

Next, we recall the definition of the \dfn{Bruhat--Chevalley order} on $\W$. For elements $u,w \in \W$, we write $u \xrightarrow{t} w$ if $u^{-1}w = t \in T$ and $\ell(u) < \ell(w)$. We define $u \leq_{BC} w$ if there exists a chain 
\[u = u_0 \rightarrow u_1 \rightarrow \cdots \rightarrow u_k = w\] 
of such relations. This chain description of the Bruhat--Chevalley order is not usually practical for direct computation. However, the following theorem gives us an equivalent characterization in terms of reduced words, which will be useful in later sections.
 \begin{theorem}[Subexpression Property,~{\cite[Theorem 2.2.2]{bjornerBrenti}}]\label{thm:bruhat-subexp}
Let $w=s'_{1}s'_{2}\cdots s'_{q}$ be a reduced word. Then $u \leq_{BC} w$ in the Bruhat--Chevalley order if and only if there is a reduced word \[u=s'_{i_1}s'_{i_2}\cdots s'_{i_k}, \quad 1\leq i_1<i_2 <\ldots < i_k\leq q,\]
obtained by deleting letters from $s'_1s'_2\cdots s'_q$.
 \end{theorem}
Consequently, $u\BC w$ if and only if \emph{some} (equivalently, \emph{every}) reduced word for $w$ contains a reduced expression for $u$ as a subexpression.

For $J\subseteq S$, let $\W_J$ be the subgroup of $\W$ generated by $J$. Subgroups of this form are called \dfn{parabolic}. Let
\[\W^J=\set{w\in\W \mid w \bc ws \text{ for all } s\in J}.\]
By \cite[Lemma 2.4.3]{bjornerBrenti}, $\W^J$ consists exactly of the elements $w\in\W$ for which no reduced expression of $w$ ends with a letter from $J$. By \cite[Corollary 2.4.5(i)]{bjornerBrenti}, these are also the unique shortest-length representatives of the left cosets $w\W_J$. Due to this bijection with $\W/\W_J$, the set $\W^J$ is abusively called a \dfn{(right) quotient}.

\begin{prop}
\label{prop:parabolic_decomp}
    Let $J \sus K \sus S$. Then the following hold:
    \begin{enumerate}
        \item  Every $w\in \W^J$ has a unique factorization $w=w^K\cdot w_K$ such that $w^K\in \W^K$ and $w_K \in \W_K \cap \W^J$. 
        \item For this factorization, $\ell(w)=\ell(w^K)+\ell(w_K)$.
    \end{enumerate}
\end{prop}
\begin{proof}
By~\cite[Proposition 2.4.4]{bjornerBrenti}, every $w \in \W$ factors uniquely as $w = w^K \cdot w_K$ with $w^K \in \W^K$, $w_K \in \W_K$, and $\ell(w) = \ell(w^K) + \ell(w_K)$. Thus, it remains only to show that if $w \in \W^J$, then $w_K \in \W^J$. 

Suppose to the contrary that $w_K \notin \W^J$. Then some reduced expression for $w_K$ ends in a generator $s \in J$. Because $\ell(w) = \ell(w^K) + \ell(w_K)$, concatenating a reduced expression for $w^K$ with this reduced expression for $w_K$ yields a reduced expression for $w$ that also ends in $s \in J$, contradicting that $w \in \W^J$.
\end{proof}

\begin{remark}
    The decomposition described by the previous proposition is called a \dfn{(right) parabolic decomposition (with respect to $K$)}. Note that in the literature (excepting \cite{richmondSlofstra}), this definition usually only considers the case with $J=\es$.
\end{remark}
\begin{defn}
Let $S(w):=\set{s\in S \mid s\BC w}.$  We will say that the parabolic decomposition of $w\in \W^J$ with respect to $K$ is \dfn{Grassmannian} if $\ab{K \cap S(w)}=\ab{S(w)}-1$ and $K$ is a maximal proper subset of $S$, that is $\ab{S-K}=1$.\footnote{Note that unless $w$ is trivial, $\ab{S(w)-J}>0$ since $w$ cannot end in any $s\in J$. In \cite{richmondSlofstra}, the authors only insist on $\ab{K \cap S(w)}=\ab{S(w)}-1$, so that $K$ is maximally proper in $S(w)\cup K$.} Likewise, any maximal proper subset $K\subset S$ will be called Grassmannian. 
\end{defn}

Note that we will apply the definition above to any finite crystallographic $(\W, S)$.

For the rest of this subsection, we fix a Coxeter group $\W$ of type $A_{n-1}$ so $\W\cong \mc{S}_n$ is generated by simple reflections $S=\set{s_1,\dots, s_{n-1}}$. We interpret $s_i$ as the \dfn{transposition} $(i, i+1)$ in the action of $\mc{S}_n$ on $[n]=\set{1,\dots,n}$. 

In the one-line notation for permutations $w=w_1w_2 \ldots w_n \in \mc{S}_n$, the right multiplication by a generator $w\cdot s_i$ swaps the entries in positions $i$ and $i+1$. On the other hand, the left multiplication $s_i\cdot w$ swaps the positions of values $i$ and $i+1$ in one-line notation. 

The \dfn{inversion set} of $w$ is defined as \[\operatorname{inv}(w) = \set{(i,j) \;\mid \; i<j,\; w_i >w_j}.\] Then $\ell(w) = |\operatorname{inv}(w)|$ and $D_R(w) = \set{s_i \in S \;\mid \; w_i > w_{i+1}}$.

We fix $1\leq k<n$ and $J=\{s_{k+1},\dots,s_{n-1}\}$, so that $\W_J\cong\mc{S}_{n-k}$ and $\W^J$ has size $(n)_k:=\frac{n!}{(n-k)!}=k!\binom{n}{k}$, recalling the falling factorial notation. 

We now recall the one-line notation of elements in $\W^J$.
\begin{lemma} \label{lem:min_reps}
    The set of shortest-length coset representatives $\W^J$ for $\W/\W_J$ consists of permutations $w=w_1\dots w_n$ such that $w_{k+1}<w_{k+2}<\dots<w_n$.
\end{lemma}
\begin{proof}
  Recall that $w\bc ws$ for all $w\in \W^J$ and $s\in J$.
     As the one-line notation of $ws_i$ is obtained by swapping $w_i$ and $w_{i+1}$ in the one-line notation of $w$, $ws_i>w$ and hence $\ell(ws_i)=\ell(w)+1$ if and only if $w_i< w_{i+1}$. Thus,
     \[ w<_{BC}ws_i \text{ for all } k+1\leq i \leq n-1 \iff w_{k+1}<w_{k+2}<\dots <w_n,\]
     completing the proof.
\end{proof}

\subsection{BP decompositions}
Here we describe a particular instance of parabolic decompositions that will be applied in \Cref{sec:results.BP}. In the literature these decompositions are often called \emph{Billey--Postnikov decompositions}, terminology which was implicitly introduced in \cite{ohYoo} for a special case applied in \cite{billeyPostnikov} to study the (rational) smoothness of Schubert varieties. The terminology was made explicit and its notions were further developed by Richmond--Slofstra in \cite{richmondSlofstra}, which we follow closely. Nevertheless, the ``BP decompositions" in question might just as well be called \emph{bundled parabolic decompositions} as we now recall. We refer the reader to \cite{borelLAG, humphreysLAG, springerLAG} for background on algebraic groups and to \cite{billeyLakshmibai} for Schubert varieties.

A finite crystallographic Coxeter system $(\W, S)$ can be realized as the Weyl group $N_G(T)/T$ of a complex connected reductive algebraic group $G$, where $T$ is a chosen maximal torus. Here, $N_G(T)$ denotes the normalizer of $T$ as a subgroup of $G$, and a choice of Borel subgroup $B\supseteq T$ determines a set of simple roots corresponding to the generators $S$. In this setting, any subset $I\sus S$ determines a parabolic subgroup $P_I\sus G$ containing $B$. The $B$-orbits of any partial flag variety $G/P_I$ are called \dfn{Schubert cells}, and their closures 
\[X^I(w):=\ol{BwP_I/P_I}\]
are called \dfn{Schubert varieties}. We will say the Schubert variety $X^I(w)$ is \emph{Grassmannian} if $I\subset S$ is a maximal proper subset. One consequence of the parabolic decomposition of \Cref{prop:parabolic_decomp} is that for any $w=w^Iw_I\in \W$, we have $X^I(w)=X^I(w^I)$. This is because $w_I$ can be realized as an element of $\W_I$, the Weyl group of $P_I$.

Given subsets $J\sus K\sus S$ with parabolic subgroups $P_J\sus P_K\sus G$, the associated partial flag varieties come with a $B$-equivariant natural projection map 
\[\pi:G/P_J \to G/P_K.\]
The restriction of this map to Schubert varieties gives a projection
\[\pi_w: X^J(w)=X^J(w^J)\to X^K(w)=X^K(w^K).\]
The following can be taken as a definition of BP decomposition, from \cite[Theorem 3.3]{richmondSlofstra}.
\begin{defn}\label{def:BP}
    Let $(\W, S)$ be the Coxeter system realized as the Weyl group determined by the triple $(G,B,T)$. Further let $w\in \W^J$ for $J\sus K\sus S$ with unique parabolic decomposition $w=vu$, $v\in \W^K$, $u\in \W_K\cap \W^J$. If the natural projection map 
    \[\pi_w: X^J(w)\to X^K(v)\]
    has the structure of a fiber bundle with fibers $X^J(u)$, then we say $w=vu$ is a \dfn{BP decomposition} with respect to $(J,K)$.
\end{defn}
In other words, a parabolic decomposition is BP if and only if the corresponding projection map is a locally trivialized fiber bundle map. We will say that a BP decomposition is \emph{Grassmannian} if the underlying parabolic decomposition is Grassmannian. The main consequences of this definition are stated in the following consolidated theorem of Richmond--Slofstra. To state it, we will denote intervals in Bruhat--Chevalley order as $$[x,y]:=\set{w\in\W \mid x\BC w \BC y}$$ and $$[x,y]^I:=[x,y]\cap \W^I$$
for arbitrary $x,y\in \W$. Further let the \dfn{Poincar\'e polynomial} of $X^I(w)$ be defined as
\[P_w^I(q):= \sum_{y \in [e,w]^I} q^{\ell(y)},\] that is the rank generating polynomial of the lower interval $[e,w]^I$. 
\begin{theorem} \label{thm:RS_BP}
Let $(\W, S)$ be the Weyl group for $(G,B,T)$ with associated Schubert varieties $X^I(w)$.  Then the following are equivalent for $J\sus K \sus S$:
\begin{enumerate}
    \item $w=vu$ is a BP-decomposition with respect to $(J,K)$.
    \item The natural projection map $\pi_w:X^J(w)\to X^K(v)$ is equidimensional.
    \item There is an isomorphism in singular cohomology with complex coefficients,
    \[H^*(X^J(w))\cong H^*(X^K(v)) \ox H^*(X^J(u))\] as $H^*(X^K(v))$-modules.
    \item There is a factorization of Poincar\'e polynomials,
    \[P^J_w(q)=P^K_v(q)\cdot P^J_u(q).\]
    \item The multiplication map $(x,y)\mapsto x\cdot y$ from the product
    \[[e,v]^K \x [e,u]^J \to [e,w]^J\]
    defines an order-preserving bijection. 
\end{enumerate}
\end{theorem}
\begin{proof}
    The equivalence (1) $\iff$ (2) is \cite[Proposition 4.7]{richmondSlofstra}. Item (3) is a consequence of our \Cref{def:BP} and the Leray--Hirsch theorem; its equivalence is proven in \cite[Corollary 3.5]{richmondSlofstra}. Item (4) is the definition of BP decomposition given by Richmond--Slofstra, and equivalence with the fiber bundle condition is their Theorem 3.3. Finally, equivalence with item (5) follows from their Proposition 4.2(b)---while that proposition statement only mentions surjectivity, in the proof they remark that the map is always injective. That it is order-preserving follows from the subexpression characterization of Bruhat--Chevalley order.
\end{proof}
\begin{remark}
    In \cite{richmondSlofstra}, the second factor in the domain of the map in item (5) of the previous theorem is taken as 
    \[[e,u]^J_K:=[e,u]\cap \W^J\cap \W_K.\]
    However since $u\in \W_K$, the Subexpression Property (\Cref{thm:bruhat-subexp}) easily implies these intervals are the same. Less obviously, by \cite[Lemma 4.8]{richmondSlofstra} the Schubert varieties associated to $u$ as an element of either $\W^J$ or $(\W_K)^J$ are also isomorphic.
\end{remark}

 The BP decomposition provides a useful framework for studying (rational) smoothness of Schubert varieties since $X^J(w)$ for $w$ with BP decomposition $w=vu$ is (rationally) smooth if and only if $X^K(v)$ and $X^J(u)$ are also (rationally) smooth \cite[Theorem 3.3(1)]{richmondSlofstra}. Recall that thanks to work of Carrell--Peterson \cite{carrellPeterson}, a Schubert variety $X^I(w)$ can be characterized as \dfn{rationally smooth} if and only if the interval $[e,w]^I$ is rank-symmetric.\footnote{For background on the geometric notion of rational smoothness and its relationship with Kazhdan--Lusztig polynomials, see \cite[Chapter 6]{billeyLakshmibai} or \cite[Section 12.2]{kumar}.} Equivalently, the Poincar\'e polynomial $P^I_w(q)$ is palindromic. We conclude this subsection by stating one further result we will need in \Cref{sec:results.BP}.

\begin{theorem}[{\cite[Theorem 3.6]{richmondSlofstra}}]
\label{thm:RS_Grassmann}
    Let $w\in \W^J$ for finite, crystallographic $\W$. If $X^J(w)$ is rationally smooth, then $w$ has a Grassmannian BP decomposition for some maximal proper $K\subset S$ containing $J$.
\end{theorem}

\subsection{Codes of permutations} 
\label{sec:background.codes}
The notion of the \emph{code} of a permutation in $\mc{S}_n$ was implicitly generalized to the setting of classical Weyl groups in \cite{stanley1980} to argue for the existence of SCDs in those types. Conceiving codes as order-preserving bijections from a product of chains to the poset $(\W, \BC)$ for arbitrary Coxeter groups, Bolognini--Sentinelli \cite{bologniniSentinelli2025} explicitly described these codes in classical type and analyzed when the map from chain products is a bijection on lower intervals. Generalized codes for lower intervals $[e,w]$ of rationally smooth elements $w\in \W$ were also treated by Gaetz--Gao \cite{gaetzGao} in connection with (complete) BP decompositions. 
We will comment further on this perspective at the end of this subsection and in \Cref{sec:final.codes}. First, we recall the classical notion of the code of a permutation and related constructions that will be necessary in \Cref{sec:results.A}. 

For the rest of \Cref{sec:background}, unless otherwise specified, we again fix a Coxeter group $\W$ of type $A_{n-1}$ so $\W\cong \mc{S}_n$ with simple reflections $S=\set{s_1,\dots, s_{n-1}}$. We first recall the construction of \emph{codes} of permutations, also known as \emph{inversion tables} \cite[p.~36]{stanleyEC1}. 

If $w=w_1\dots w_n$ is the one-line notation of a permutation $w\in \mc{S}_n$, then the \dfn{Rothe diagram} $\D(w)$ can be defined as 
\[\D(w)=\set{(i,w_j) \mid i<j\text{ and }w_i>w_j}.\]
 The Rothe diagram can be depicted as a collection of boxes constructed from the permutation matrices associated to elements of $\mc{S}_n$ as described below.\footnote{Note that conventions in the literature on Rothe diagrams vary, sometimes preferring row permutation matrices and collections of boxes given in plane coordinates.}

For $w\in \mc{S}_n$, let $M_w=(a_{ij})$ be the (column) permutation matrix, defined by 
\[ a_{ij}=\begin{cases}
   1 & \text{ if } w_j=i,\\
   0 & \text{ otherwise.}
\end{cases}
    \]
From this matrix, at each entry which is equal to 1, one extends ``death rays" to all of the positions directly south or east of a given 1. The positions of the matrix which are not touched by any death rays constitute the boxes of the Rothe diagram. To distinguish the visual representation just described from the set-theoretic description of the Rothe diagram above, we will call it the \dfn{ray diagram} of the permutation. Observe that the pairs $(i,w_j)$ of the Rothe diagram index the matrix coordinates $a_{w_j i}$ of the surviving boxes in the ray diagram.

\begin{example}\label{ex:rothe-diagram}
Consider the permutation $w=6132457$. Its corresponding ray diagram is depicted in Figure~\ref{fig:ray-diagram}.
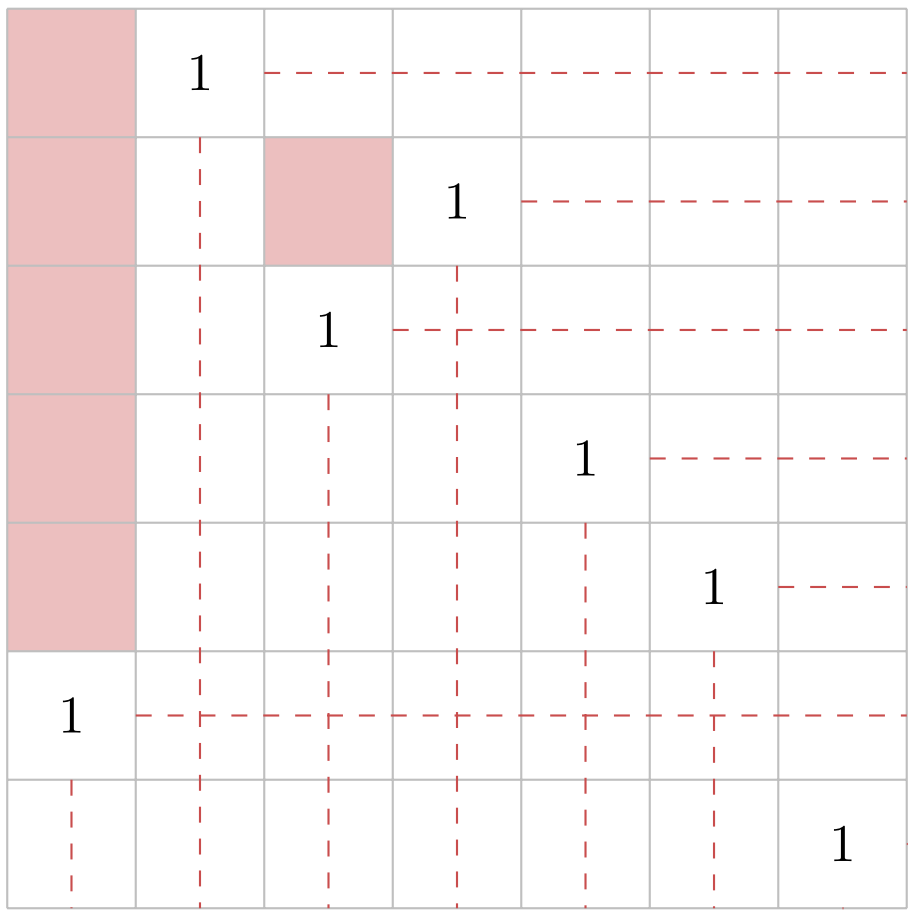
\begin{figure}[htb!]
        \centering
        \begin{tikzpicture}[scale=0.85]
            \foreach \r/\c in {1/1,2/1,3/1,4/1,5/1,2/3}{
            \fill[darkred!25] (\c-1, {7-\r}) rectangle (\c, {8-\r}); 
            }
            \draw[step=1cm, gray!50] (0,0) grid (7,7);
             \foreach \r/\c in {6/1,1/2,3/3,2/4,4/5,5/6,7/7}{
            \draw[dashed, darkred!70] (\c-0.5, {7-\r}) -- (\c-0.5, 0); 
            \draw[dashed, darkred!70] (\c, {7.5-\r}) -- (7, {7.5-\r}); 
            }
            \foreach \r/\c in {6/1,1/2,3/3,2/4,4/5,5/6,7/7}{
            \node at (\c-0.5,{7.5-\r}) {$1$};
            }
        \end{tikzpicture}
        \caption{The ray diagram of $w=6132457$. The dashed lines are the death rays extending south and east from each $1$. The shaded boxes are the surviving cells, giving the Rothe diagram of $w$.}
        \label{fig:ray-diagram}
    \end{figure}
    Thus, the Rothe diagram is $\set{(1,1),(1,2), (1,3),(1,4),(1,5),(3,2)}$.
\end{example}

\noindent
The \dfn{code} of a permutation $w\in \mc{S}_n$ is defined as the sequence 
\[\code(w)=(a_1,\dots,a_n)\]
where $a_i$ is the number of boxes in the $i^{th}$ column of the ray diagram of $w$. Alternatively, $a_i$ counts the inversions $(i,j)\in \inv(w)$ with fixed $i$ and some $j$; that is, $$a_i=|\{j \in \{i+1,\dots,n\}\mid w_j<w_i\}|.$$

\begin{example}
The code of the permutation $w=6132457$ from Example~\ref{ex:rothe-diagram} is $\code(w)=(5,0,1,0,0,0,0)$.
\end{example}

The code of a permutation is also useful for constructing a reduced word in terms of the Coxeter generators, as we now recall. For any positive integer $j$, let $\ul{j}$ be the chain poset which consists of the usual ordering of the natural numbers $\set{0,1,\dots, j-1}$, where the rank of $i\in \ul{j}$ is $i$. We will also use this notation to refer to the underlying set. Let $\mb{n}_k$ denote the product poset
\[\mb{n}_k=\ul{n}\x \ul{n-1}\x  \cdots \x \ul{n-k+1}. \]

Following the conventions of \cite{lambertRabelo2026}, we can regard an element $\alpha=(a_1,\dots,a_n)\in\mb{n}_n$ as a weak composition. Its \dfn{composition diagram} is the set of cells $(i,j)$ with $1\leq j\leq a_i$, drawn in English notation, so that row $i$ has $a_i$ left-justified cells. Since $a_i\leq n-i$, the composition diagram fits inside the staircase partition $\la_{st_n}=(n-1,n-2,\dots,1)$. The order relation $\preceq$ on $\mb{n}_n$ is then given by containment of composition diagrams. As codes of $\mc{S}_n$ elements are in bijection with such compositions, the composition diagram of $\alpha=\code(w)$ is obtained by reflecting the ray diagram of $w$ across the main diagonal and left-justifying its cells.

\begin{defn} \label{def:codeword}
    Let $\alpha\in\mb{n}_n$ and label the cells of the weak composition diagram inside $\lambda_{st_n}$ using matrix coordinates $(i,j)$, so that $1\leq i,j \leq n-1$ and $i+j\leq n$. The \dfn{codeword} associated to $\alpha$, denoted $\ul{w}(\alpha)=s_{c_1}\cdots s_{c_{\ell(w)}}$, is defined as follows:
\begin{enumerate}
    \item Label each box $(i,j)$  with $s_{i+j-1}$. 
    \item For each row $i$, construct a word $\ul{r}_\alpha(i)$ by selecting the shaded boxes from right to left within each row and concatenating the labels.
    \item Define $\ul{w}(\alpha)=\ul{r}_\alpha(1)\cdots \ul{r}_\alpha(n-1)$,    
    concatenating the word of each row from top to bottom. 
\end{enumerate} 
Given an element $\alpha\in \mb{n}_n$, we denote by $w(\alpha)$ the group element associated to the codeword $\ul{{w}}(\alpha)$.
\end{defn}

\begin{example}\label{ex_codeword}
Recall the permutation $w=6132457$ from Example~\ref{ex:rothe-diagram} with code $\alpha = (5,0,1,0,0,0,0)$. The shaded cells of its composition diagram are
$$\set{(1,1),(1,2), (1,3),(1,4),(1,5),(3,1)}.$$ Labeling each box according to~\Cref{def:codeword}(1) and reading these right to left within each row yields $\ul{r}_\alpha(1)=s_5 s_4 s_3 s_2 s_1$ and $\ul{r}_{\alpha}(3)=s_3$. Concatenating these subexpressions gives the codeword $\ul{w}(\alpha) = s_5 s_4 s_3 s_2 s_1 s_3.$
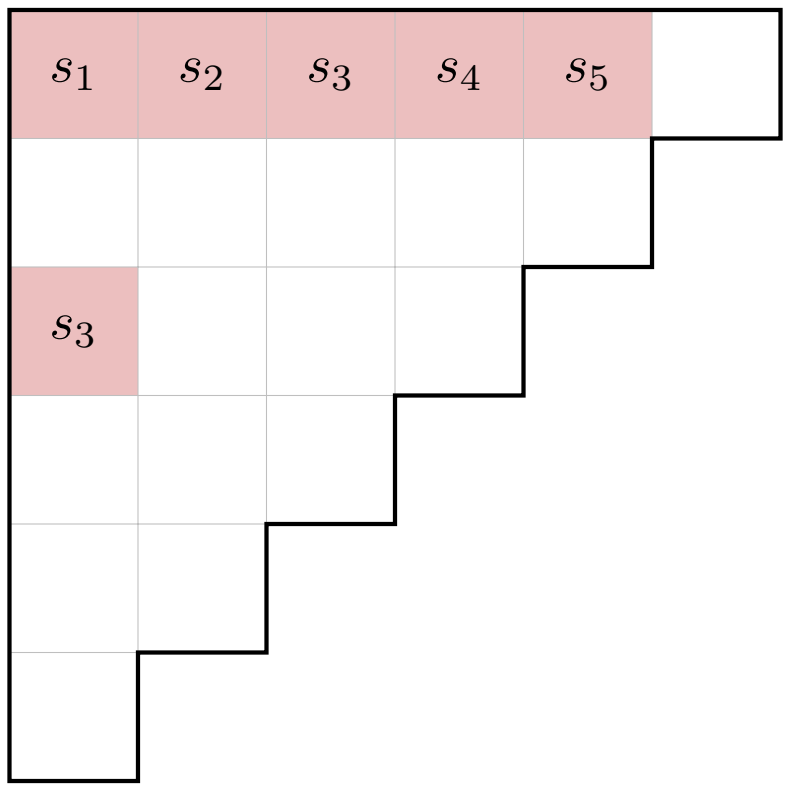
\begin{figure}[htb!]
    \centering
    \begin{tikzpicture}[scale=0.85]
        \fill[darkred!25] (0,6) rectangle (1,7);
        \fill[darkred!25] (1,6) rectangle (2,7);
        \fill[darkred!25] (2,6) rectangle (3,7);
        \fill[darkred!25] (3,6) rectangle (4,7);
        \fill[darkred!25] (4,6) rectangle (5,7);
        \fill[darkred!25] (0,4) rectangle (1,5);
        \begin{scope}
            \clip (0,7) -- (6,7) -- (6,6) -- (5,6) -- (5,5) -- (4,5) -- (4,4) -- (3,4) -- (3,3) -- (2,3) -- (2,2) -- (1,2) -- (1,1) -- (0,1) -- cycle;
            \draw[gray!50, very thin] (0,0) grid (6,7);
        \end{scope}
        \draw[thick] (0,7) -- (6,7) -- (6,6) -- (5,6) -- (5,5) -- (4,5) -- (4,4) -- (3,4) -- (3,3) -- (2,3) -- (2,2) -- (1,2) -- (1,1) -- (0,1) -- cycle;
        \node at (0.5,6.5) {$s_1$};
        \node at (1.5,6.5) {$s_2$};
        \node at (2.5,6.5) {$s_3$};
        \node at (3.5,6.5) {$s_4$};
        \node at (4.5,6.5) {$s_5$};
        \node at (0.5,4.5) {$s_3$};
    \end{tikzpicture}
    \caption{The staircase partition $\lambda_{st_7} = (6,5,4,3,2,1)$, with the codeword $\ul{w}(\alpha) = s_5 s_4 s_3 s_2 s_1 s_3.$} 
    \label{fig:staircase-codeword}
\end{figure}
\end{example}

\begin{lemma}[{Cf. \cite[Remark 2.1.9]{manivel}}] \label{lem:staircase bijection}
The code of a permutation defines a bijection
\[\code: \mc{S}_n \to \ul{n}\x \ul{n-1} \x \cdots \x \ul{1}=\mb{n}_n,\]
whose inverse is given by $\code^{-1}(\alpha)=w(\alpha)$ where $w(\alpha)$ is as in \Cref{def:codeword}. 
\end{lemma}
\begin{proof}
Let $\code(w)=(a_1,\dots,a_n)$ for an arbitrary element $w\in \mc{S}_n$. By definition, $a_i$ is the number of inversions of the form $(i,k)\in \inv(w)$ with fixed first coordinate $i$. Because there are $n-i$ positions following $i$, there are at most $n-i$ such inversions, ensuring $a_i \in \ul{n-i+1}$ and $\code(w) \in \mb{n}_n$. Since $\mc{S}_n$ and $\mb{n}_n$ are equinumerous, it suffices to show that $\code$ is injective. Following~\cite[Proposition~2.1.2]{manivel}, injectivity follows immediately from the recursive reconstruction algorithm for $w$ from $\code(w)$: for each $i$, set $w_i$ to be the $(a_i+1)$-th smallest element of the sequence obtained from $1,\dots,n$ by omitting $w_1,\dots, w_{i-1}$.

To verify that the explicit codeword construction of \Cref{def:codeword} yields the inverse map $\code^{-1}(\alpha)=w(\alpha)$, we show that $\code(w(\alpha))=\alpha$ for any $\alpha=(a_1,\dots,a_n) \in \mb{n}_{n}$. 

For a given row $i$ with $a_i > 0$, setting $j = i + a_i$ we have the word
\[
\ul{r}_\alpha(i) = s_{j-1} s_{j-2} \cdots s_i,
\]
which corresponds to the cycle $\sigma$ with $\sigma_i = j$, $\sigma_{i+1} = i$, and $\sigma_k = k-1$ for $i+2 \leq k \leq j$. This cycle has precisely $a_i = j - i$ inversions which are all of the form $(i,k)$ for $i+1\leq k \leq j$.  

 Let $w^{(i)}$ denote the product of the generators corresponding to the word $\ul{r}_\alpha(1)\cdots \ul{r}_\alpha(i)$, and let $$\code(w^{(i)})=(a_1^{(i)},\dots,a_n^{(i)}).$$ 
 We argue by induction on $i$ that 
 $$(a_1^{(i)},\dots a_i^{(i)}, a_{i+1}^{(i)},\dots, a_n^{(i)})=(a_1,\dots, a_i, 0,\dots ,0)$$ for all $0\leq i\leq n-1$, so that $\code(w^{(n-1)})=\code(w(\alpha))=\alpha$. For the base case $i=0$, the statement holds trivially since $w^{(0)}$ is the identity permutation, whose code is $(0,\dots,0)$.
 
 Now suppose the statement holds for some $i>0$. First we argue that $a_s=a_s^{(i+1)}$ for all $s\leq i$. Observe that $\ul{r}_\alpha(i+1)$ consists exclusively of generators $s_m$ with $m \geq i+1$. Right-multiplication by these generators leaves the entries of the one-line notation of $w^{(i)}$ at position $s\leq i$ fixed and merely permutes the entries in positions $\{i+1,\dots,n\}$. Consequently, the number of entries to the right of position $s$ that are smaller than $w^{(i)}_s=w^{(i+1)}_s$ remains equal to $a^{(i)}_s=a_s$ in $\code(w^{(i+1)})$.
 
Next we argue that $a^{(i+1)}_{i+1}=a_{i+1}$ and that $a^{(i+1)}_t=0$ for $t>i+1$. Since $a^{(i)}_{t}=0$ for all $t>i$, the sequence of symbols $(w^{(i)}_{i+1}, \dots,w^{(i)}_{n})$ is increasing. If $\ul{r}_\alpha(i+1)$ is empty, then we are done. Otherwise, multiplication on the right by the cycle corresponding to $\ul{r}_\alpha(i+1)$ acts by sending $w^{(i)}_j$ to position $i+1$ for some $j>i+1$, and maintains the relative order of all other symbols $w^{(i)}_t$, $t\geq i+1$. Thus, it introduces exactly $a_{i+1}$ inversions of the form $(i+1,k)$ and leaves the sequence of symbols $(w^{(i+1)}_{i+2},\dots w^{(i+1)}_{n})$ increasing. This completes the proof.
\end{proof}

\begin{example}
Consider $\alpha = (2,0,1,0,0) \in \mb{5}_5$ whose nonzero entries are $a_1=2$ and $a_3 =1$. Using the construction of $\ul{{w}}(\alpha)$ from \Cref{def:codeword}, we have
$$\ul{r}_\alpha(1) = s_2\,s_1, \quad \ul{r}_\alpha(3) = s_3, \quad  \text{and} \quad\ul{{w}}(\alpha) = \ul{r}_\alpha(1) \cdot \ul{r}_\alpha(3) = s_2\,s_1\,s_3.$$

Applying this to the identity yields $w(\alpha) = 31425$. Moreover, since $w_1=3$ produces two inversions and $w_3=4$ produces one inversion, $\code(31425) = (2,0,1,0,0) = \alpha$, as in Lemma~\ref{lem:staircase bijection}. 
\end{example}
\begin{remark} \label{rk:code.BPseq}
     The code of a permutation in $\mc{S}_n$ can be understood as an artifact of repeated application of (bundled) parabolic decomposition as we now describe. For $0\leq i\leq n-1$ let $J_i=\set{s_{n-i},\dots, s_{n-1}}$ so that $J_0=\varnothing$ and $J_{n-1}=S$. Then
    \[\W_{J_i}^{J_{i-1}}=\set{e, s_{n-i},s_{n-i+1}s_{n-i},\dots, s_{n-1}\cdots s_{n-i+1}s_{n-i}},\]
    so we see that the Bruhat--Chevalley order on this poset is a chain of order $i+1$. As in \cite[Corollary 2.4.6]{bjornerBrenti}, there is then an order- and rank-preserving bijection given by the multiplication map
    \begin{align*}  
        \W_{J_{n-1}}^{J_{n-2}}\x\W_{J_{n-2}}^{J_{n-3}} \x \cdots \x \W_{J_1}^{J_0}& \to \W \\
        (v_{n-1},\dots, v_1)&\mapsto v_{n-1}\cdots v_1.
    \end{align*}
Letting $w=v_{n-1}\cdots v_1$ with $\code(w)=\alpha=(a_1,\dots,a_n)$, we then see that $a_i=\ell({v_{n-i}})$ and $w^{(i)}=v_{n-1}\cdots v_{n-i}$ in the notation of the proof of \Cref{lem:staircase bijection}. Furthermore, supposing that each $v_{i}$ is the longest element of $\W^{J_{i}}_{J_{i-1}}$, then $v_{i}(v_{i-1}\cdots v_1)$ is a Grassmannian BP decomposition of the longest element of $\W_{J_{i}}$ for each $i$. 
\end{remark}

\begin{remark}   \label{rk:code.restrict}
      Let $\hat{0}:=(0,\dots,0)$ denote the minimal element of $\mb{n}_n$ so that $[\hat{0},\code(w)]$ is the lower interval of $\mb{n}_n$ for some composition $\code(w)$ associated to $w\in \mc{S}_n$. While it is known \cite{gasharov, lascoux, billey98} that intervals $[e,w]$ for (rationally) smooth permutations $w$ admit order-preserving bijections from products of chains (that is, the g-codes of the introduction), it is not the case that such bijections always come from $[\hat{0},\code(w)]$. 
      
      To see this, consider $w=312$ with $\code(w)=(2,0,0)$ so that $[\hat{0},\code(w)]\cong \ul{3}$, but $[e,w]\cong \ul{2}\x\ul{2}$. On the other hand $w^{-1}=231$ has $\code(w^{-1})=(1,1,0)$ with $[e,w^{-1}]\cong [\hat{0}, \code(w^{-1})]\cong\ul{2}\x\ul{2}$.    
      The set of smooth elements such that the code for $\mc{S}_n$ yields a g-code for $[e,w]$ is described in \cite[Proposition 7.8]{bologniniSentinelli2025}.
\end{remark}

\subsection{Rank-control matrices} Among the many criteria available to determine Bruhat--Chevalley order relations for permutations, we describe here one that can be easily related to $\code(w)$ due to its construction from a permutation matrix. Given $w\in \mc{S}_n$, the (northwest) \dfn{rank-control matrix} associated to $w$ is an $n\x n$ matrix denoted $r_w$ with $ij$-entry 
\[r_w(i,j)=\ab{\set{(w_k,k) \mid w_k\leq i,\ k\leq j }}.\]

\begin{lemma}[{\cite[Proposition~2.1.12]{manivel}}] \label{lem:rk_bc}
    For $v,w \in \mc{S}_n$, $v\BC w$ if and only if for all $1\leq i,j \leq n$ we have 
    \[r_v(i,j)\geq r_w(i,j).\]
\end{lemma}
We write $r_v\geq r_w$ to denote that $r_v(i,j)\geq r_w(i,j)$ for all $i,j$. That is, 
\[v\BC w \iff r_v\geq r_w.\]
Thinking of the associated permutation matrix as a non-attacking rook placement on a square board, the rank control matrix also counts the number of rooks weakly northwest of a given position. Note that our convention again differs from that of \cite{manivel} by preferring column permutation matrices. However, the northwest condition is unaffected by this modification. See also \cite[Theorem 2.1.5]{bjornerBrenti}.

\begin{remark}\label{rk:tableau}
The criterion in~\Cref{lem:rk_bc} can be rephrased in terms of tableaux.  Following the convention of~\cite[Theorem~2.6.3]{bjornerBrenti}, for $w\in\mc{S}_n$ and $1\leq k\leq n$ let $w_{1,k}<\dots<w_{k,k}$ denote the increasing rearrangement of $w_1,\dots,w_k$. The tableau of $w$ is the array whose $k$th row is $(w_{1,k},\dots,w_{k,k})$. Note that $r_w(i,j)$ is the number of entries of $\set{w_{1,j}, \dots, w_{j,j}}$ that are at most $i$. Hence for $v \in\mc{S}_n$, $r_v \geq r_w$ if and only if for each $j$, the $\ell$th smallest entry of $\set{v_{1}, \dots, v_{j}}$ is at most the $\ell$th smallest entry of $\set{w_{1}, \dots, w_{j}}$. Thus, $v \BC w$ if and only if $v_{i,k} \leq w_{i,k}$ for all $1\leq i \leq k \leq n$, that is, if and only if the tableau of $v$ is entrywise at most the tableau of $w$. Moreover, to check whether $v \BC w$ it is enough to compare the rows $k$ with $s_k \in D_R(v)$. For another tableau criterion, see, for instance,~\cite[Theorem 5A, Proposition~7.3]{proctor82}.

A tableau criterion for $(\mc{S}_n,\BC)$ was first defined by Ehresmann, and  Deodhar extended this perspective to all Coxeter groups. 
\end{remark}

\begin{example}
Consider the permutations $w=6132457$ and $v=3162457$ in $\mc{S}_7$. 
Applying the definition of the rank-control matrix to all entries, we get: 
\[
r_w = \begin{pmatrix}
0 & 1 & 1 & 1 & 1 & 1 & 1\\
0 & 1 & 1 & 2 & 2 & 2 & 2\\
0 & 1 & 2 & 3 & 3 & 3 & 3\\
0 & 1 & 2 & 3 & 4 & 4 & 4\\
0 & 1 & 2 & 3 & 4 & 5 & 5\\
1 & 2 & 3 & 4 & 5 & 6 & 6\\
1 & 2 & 3 & 4 & 5 & 6 & 7
\end{pmatrix}
\qquad \text{and} \qquad
r_v = \begin{pmatrix}
0 & 1 & 1 & 1 & 1 & 1 & 1\\
0 & 1 & 1 & 2 & 2 & 2 & 2\\
1 & 2 & 2 & 3 & 3 & 3 & 3\\
1 & 2 & 2 & 3 & 4 & 4 & 4\\
1 & 2 & 2 & 3 & 4 & 5 & 5\\
1 & 2 & 3 & 4 & 5 & 6 & 6\\
1 & 2 & 3 & 4 & 5 & 6 & 7
\end{pmatrix},
\]
respectively. Since $r_v(i,j)\geq r_w(i,j)$ for all $i,j$, $r_v\geq r_w$. Hence, $v \BC w$ by Lemma~\ref{lem:rk_bc}.

On the other hand, since $v_1 > v_2$ and $v_3 > v_4$ are the only descents of $v$, we have $D_R(v)=\set{s_1, s_3}$. Hence, by~\Cref{rk:tableau}, it is enough to compare the rows $k=1$ and $k=3$ of the tableaux of $v$ and $w$. The tableaux of these rows, from top to bottom, are
\[
\begin{ytableau}
6\\
1 & 3 & 6
\end{ytableau}
\qquad\text{and}\qquad
\begin{ytableau}
3\\
1 & 3 & 6
\end{ytableau}
\]
for $w$ and $v$, respectively. From the tableaux, we have $v_{1,1} =3 \leq w_{1,1}=6$ and for $k=3$ we have  $(v_{1,3}, v_{2,3}, v_{3,3}) = (1,3,6) = (w_{1,3}, w_{2,3}, w_{3,3})$. Hence, $v \BC w$.
\end{example}

\begin{remark}\label{rk:code_RCM}
 Notice that the code component $a_i = w_i - r_{w}(w_i,i)$ for all $1\leq i \leq n$.  
\end{remark}

\subsection{Symmetric chain decompositions}

Here we recall some basic facts about symmetric chain decompositions (SCDs). For further background, see \cite{anderson, engel}. Given a finite poset $\mc{P}$, a \dfn{chain} $C=(c_1,\dots,c_r)$ is a tuple of elements $c_i\in\mc{P}$ such that $c_i < c_{i+1}$ for each $1\leq i \leq r-1$. A \dfn{saturated chain} is a chain $C=(c_1,\dots,c_r)$ such that $c_{i+1}$ covers $c_i$ for each $1 \leq i \leq r-1$. The \dfn{length} of a chain $(c_1,\dots,c_r)$ is $r-1$ (for a saturated chain, this is the number of covering relations it comprises).

Suppose that $\mc{P}$ is a poset with a unique minimal element $\hat{0}$ such that for every $c \in \mc{P}$, all saturated chains from $\hat{0}$ to $c$ have the same length. In this case, we call $\mc{P}$ a \dfn{ranked poset}. We define the \dfn{rank function} $\rank: \mc{P} \to \N$ by setting $\rank(c)$ to be the length of any saturated chain from $\hat{0}$ to $c$, and we let $\rank(\mc{P}) = \max_{c \in \mc{P}} \rank(c)$ denote the rank of the poset.

Given a ranked poset $\mc{P}$, a saturated chain $C=(c_1,\dots,c_r)$ is said to be a \dfn{symmetric chain} if 
\[
\rank(c_1) + \rank(c_r) = \rank(\mc{P}).
\]
A \dfn{symmetric chain decomposition} (SCD) of $\mc{P}$ is a partition of $\mc{P}$ into pairwise disjoint symmetric chains. A poset that admits an SCD is called a \dfn{symmetric chain order} (SCO).

The following lemma is one of the basic techniques for constructing SCDs of some poset from existing SCDs. A proof appears in \cite[Theorem 3.6.1]{anderson} and \cite[Theorem 5.1.5]{engel}.

\begin{lemma} \label{lem:product_SCD}
    If posets $\mc{P}$ and $\mc{Q}$ are symmetric chain orders, then so is $\mc{P} \x \mc{Q}$. 
\end{lemma}
Consequently, any product of finite chains is a symmetric chain order. The proof is based on an elementary decomposition of the product of two chains $\ul{a} \x \ul{b}$. In fact, there are two natural options for how to specify an SCD of the product, according to whether the chains ascend in the $\ul{a}$ or the $\ul{b}$ component first. An example is depicted below in Figure~\ref{fig:SCD_a} and Figure~\ref{fig:SCD_b}. We will combine \Cref{lem:product_SCD} with BP decompositions to construct SCDs of certain Bruhat--Chevalley posets.

\section{Main Results}
\label{sec:results}

We begin by showing a statement for Weyl groups which implies that in type $A$, \Cref{prob:gen} reduces to \Cref{prob:rect}. Then we proceed to describe SCDs on two Bruhat--Chevalley posets relevant to \Cref{prob:rect}. Both of these posets contain the poset $L(k,n-k)$ of partitions in a $k\x (n-k)$ rectangle as a subposet, and the latter is of the same rank. 

\subsection{SCDs and (complete) BP decompositions}
\label{sec:results.BP}
Throughout this subsection and unless indicated otherwise, $(\W,S)$ is a finite crystallographic Coxeter system; that is, the Weyl group of some complex connected reductive algebraic group.

It is obvious that a poset needs to be rank-symmetric in order to possibly admit an SCD. Our main result for this subsection suggests that, for finite Weyl groups, rank-symmetry may be a sufficient condition for any Bruhat-Chevalley interval to admit an $SCD$. The key result here is the following, which is proved by repeated application of \Cref{thm:RS_Grassmann}. A version in type $A$ was first proven for $X(w)$ by Ryan \cite{ryan} and later extended to all $X^J(w)$ by Wolper \cite{wolper}.

\begin{theorem}[{\cite[Corollary 3.7]{richmondSlofstra}}]
\label{thm:RS_RW}
        Let $w \in \W^J$ for finite crystallographic $(\W,S)$, and set $m=\ab{S(w)-J}$. Then $X^J(w)$ is rationally smooth if and only if there is a sequence of Grassmannian BP decompositions $v_i(v_{i-1} \cdots v_1)$ where: 
        \begin{enumerate}
            \item $w=v_m \cdots v_1$,
            \item $v_i \in \W_{J_i}^{J_{i-1}}$ for 
            \[ J_i= J \cup \bigcup_{j\leq i}S(v_j),\]
            \item and $w_i:=v_m\cdots v_{i+1} \in \W^{J_i}$ for $0\leq i \leq m$  (note $w_m=e$),
        \end{enumerate}
       so that there is a sequence of fiber bundle maps 
        \[X^J(w)=X^{J_0}(w_0)\to X^{J_1}(w_1)\to \cdots \to X^{J_{m-1}}(w_{m-1}) \to X^{J_m}(w_m) = {\text{pt}} \]
        where the fiber over $X^{J_{i+1}}(w_{i+1})$ is the rationally smooth Grassmannian Schubert variety $X^{J_i}(v_{i+1})$. 
\end{theorem}

A sequence of Grassmannian BP decompositions as described in the theorem above is called a \dfn{complete BP decomposition}. Our next statement unpacks the consequences of this result for an interval $[e,w]^J$. 
\begin{prop} \label{prop:RS_cor}
    With notation as in \Cref{thm:RS_RW}, $[e,w]^J$ is rank-symmetric if and only if:
    \begin{enumerate}
        \item There is an order preserving bijection from the product poset
        \[\prod_{i=1}^m [e,v_{i}]^{J_{i-1}} \to [e,w]^J,
        \]
        where each factor is a rationally smooth Grassmannian Schubert variety.
        \item There is a factorization of Poincar\'e polynomials 
        \[P_w^J(q)=\prod_{i=1}^m P_{v_i}^{J_{i-1}}(q),\]
        where the pairs $(v_i, J_{i-1})$ are as above.
        \item There is an isomorphism in singular cohomology
        \[H^*(X^J(w))\cong \bigotimes_{i=1}^m H^*(X^{J_{i-1}}(v_i))\]
        as $H^{*}(X^{J_{m-1}}(v_m))$-modules.
    \end{enumerate}
\end{prop}
\begin{proof}
    That rank-symmetry/rational smoothness implies the given statements follows directly from \Cref{thm:RS_RW} and the equivalent conditions of \Cref{thm:RS_BP}. The other direction follows from the fact that the rationally smooth Grassmannian Schubert variety factors have palindromic Poincar\'e polynomials, hence so does their product. The $H^*(X^{J_{m-1}}(v_m))$-module structure on $H^*(X^J(w))$ comes by iteratively applying the pullback on cohomology induced by each bundle map, noting that $w_{m-1}=v_m$.
\end{proof} 
\begin{theorem}\label{thm:main1}
    If the Bruhat--Chevalley poset of any rationally smooth Grassmannian Schubert variety in any irreducible Coxeter type is an SCO, then the interval $[e,w]^J$ for any $J\sus S$ and any rationally smooth Schubert variety $X^J(w)$ is an SCO. 
\end{theorem}
\begin{proof}
    Suppose $X^J(w)$ is rationally smooth. Then~\Cref{prop:RS_cor}(1) implies that there is an order-preserving bijection from a product of rationally smooth Grassmannian Schubert variety intervals $[e,v_i]^{J_{i-1}}$ to $[e,w]^J$. By the construction of~\Cref{thm:RS_RW}, we have $\ab{J_i-J_{i-1}}=1$. Let $s'_i$ denote the element of $J_i- J_{i-1}$, and let ${K_i}$ denote the connected component of the Coxeter graph of $J_i$ that contains $s'_i$. We then have $\W_{J_i}^{J_{i-1}} = \W_{K_i}^{K_i-s_i'}$ and since $\W_{K_i}$ is irreducible, the interval $[e,v_i]^{J_{i-1}}$ is isomorphic to that of a rationally smooth Grassmannian element of a irreducible Coxeter system. Thus, supposing all of these factors are SCOs, then so is $[e,w]^J$ by repeated application of \Cref{lem:product_SCD}.
\end{proof}
It is known that in type $A$~\cite{deodhar85} rational smoothness is equivalent to smoothness. Furthermore, all smooth Grassmannian Schubert varieties are sub-Grassmannian varieties \cite[Theorem 9.3.1]{billeyLakshmibai}, and so their Bruhat--Chevalley posets are of the form $L(k,n-k)$. As the rank-generating polynomials of these posets are $q$-binomial coefficients,~\Cref{prop:RS_cor}(2) implies that $P_w^J(q)$ factors as a product of $q$-binomials for all smooth $X^J(w)$. Moreover, we can deduce~\Cref{thm:intro} from the above.
\begin{cor}
    Suppose that $L(k,n-k)$ is an SCO for all $1\leq k <n$. Then the Bruhat--Chevalley order on any smooth type $A$ Schubert variety is also an SCO. 
\end{cor}
In particular, in type $A$ the order-preserving bijection of~\Cref{prop:RS_cor}(1) is of the form
\[ \prod_{i=1}^m L(k_i,n_i) \to [e,w]^J \]
for some $1\leq k_i<n_i$. From the results on SCDs of $L(k,n-k)$ mentioned in the introduction, this tells us, for instance, that if $X^J(w)$ is such that $k_i\leq 6$ for each $i$, then $[e,w]^J$ is an SCO. The known factorization of $P_w(q)$ into $q$-integers for smooth $w\in \W$ \cite{gasharov, lascoux} can then be understood as a special case where each $k_i=1$. In this case, each Grassmannian is the space of lines in a complex vector space, and hence equal to a complex projective space $\C\P^{n_i}$.

We close this subsection with the following result which applies to any finite Coxeter group.

\begin{prop}\label{prop:cox_subquotients_bij}
    Let $K \subset S$ for a finite Coxeter system $(\W,S)$ such that $\W_K=\W_I \x \W_J$ for some $I, J\subset K$. Then there is an order-preserving bijection $\W^K\x\W_{I}\to \W^J$ and hence a BP decomposition $w=vu$ where $w$ is the longest element of $\W^J$, $v$ is the longest element of $\W^K$, and $u$ the longest element of $\W_I$.
\end{prop}
\begin{proof}
    Let $I,J,K$ be as in the proposition statement. 
We define the mapping $\phi$ as follows:
\begin{equation*}\label{eq:phi_order_preserving-function}
\begin{tikzcd}[row sep=tiny,column sep=large]
\W^K\x \W_I \arrow{r}{\phi} & \W^J \\
(x,y)   \arrow[r, maps to] & \phi(x,y) := xy 
\end{tikzcd}
\end{equation*}
We must show that $\phi$ is a well-defined bijection and that it is order-preserving. 

\emph{$\phi$ is a well-defined bijection:} By the parabolic decomposition with respect to $K$ of \Cref{prop:parabolic_decomp}, any element $w \in \W^J$ can be uniquely factored as $w=xy$ where $x \in \W^K$, $y\in \W_K \cap \W^J$, and $\ell(w)=\ell(x)+\ell(y)$. 
Since $\W_K=\W_I \x \W_J$, every element $y \in \W_K$ factors uniquely as $y=zt$ where $z \in \W_I$ and $t \in \W_J$. For $v$ to be a minimal length representative modulo $\W_J$, we must have $t = e$. Therefore, $y = z \in \W_I$, and so $\W_K \cap \W^J=\W_I$. In other words, every $w \in \W^J$ factors uniquely as $w = xy$ with $x \in \W^K$, $y \in \W_I$, and $\ell(w) = \ell(x) + \ell(y)$. This guarantees $\phi$ is a bijection.

\emph{$\phi$ is order-preserving:} Let $\W^K\x \W_I$ have the product order: $(x_1, y_1) \leq (x_2, y_2)$ if and only if $x_1 \leq_{BC} x_2$ in $\W^K$ and $y_1 \leq_{BC} y_2$ in $\W_I$. Now, suppose that $(x_1, y_1) \leq (x_2, y_2)$. 

Let $x_2$ have a reduced word $s_{a_1}\cdots s_{a_p}$, and let $y_2$ have a reduced word $s_{b_1}\cdots s_{b_q}$. Because $\ell(x_2y_2)=\ell(x_2)+\ell(y_2)$, concatenating them yields a valid reduced word for the product:
\[
x_2y_2 = s_{a_1}\cdots s_{a_p} s_{b_1}\cdots s_{b_q}.
\]
It follows from \Cref{thm:bruhat-subexp} that there exists a subexpression of $s_{a_1}\cdots s_{a_p}$ evaluating to $x_1$ (since $x_1\leq_{BC} x_2$), and there exists a subexpression of $s_{b_1}\cdots s_{b_q}$ evaluating to $y_1$ (since $y_1\leq_{BC} y_2$).
Concatenating these extracted subexpressions produces a valid subexpression of the full reduced expression for $x_2y_2$ that evaluates to $x_1y_1$. Applying \Cref{thm:bruhat-subexp} once more, the existence of this subexpression confirms that $x_1y_1 \leq_{BC} x_2y_2$ in the Bruhat--Chevalley order on $\W^J$.
\end{proof}

\subsection{Permutation-labeled partitions in a rectangle}
For the remainder of this section, we return to our realization of $\mc{S}_n$ as $(\W, S)$, a Coxeter system of type $A_{n-1}$. We fix $J$ as the subset $\set{s_{k+1},\dots, s_{n-1}}$ so that $\W/\W_J\cong \mc{S}_n/\mc{S}_{n-k}$ unless indicated otherwise. The combinatorial models we now present reflect the bijection of \Cref{prop:cox_subquotients_bij}, which we will restate in \Cref{cor:order_preserving_bijection_Sn}. 
\begin{prop}\label{lemma:cosets_and_lattice_paths} \label{prop:coset_pairs}
        The coset representatives of $\W^J$ are in correspondence with pairs $(\mb{L},\pi)$ such that $\mb{L}$ is a $\Z^2$ lattice path from the origin to $(n-k,k)$ and $\pi \in \mc{S}_k$ is a permutation. 
\end{prop}
\begin{proof} The parabolic subgroup $\W_J\cong\mc{S}_{n-k}$ acts exclusively on the last $n-k$ positions. From Lemma~\ref{lem:min_reps}, $w \in \W^J$ if and only if $w_{k+1} < w_{k+2} < \dots < w_n$. Given $w \in \W^J$, let $A = \{w_1, \dots, w_k\}$ be the set of values in the first $k$ positions. We construct a bijection from $w$ to a pair $(\mb{L}, \pi)$ as follows.

\textit{The lattice path $\mb{L}$.} We encode $A$ as a path on $\Z^2$ from $(0,0)$ to $(n-k, k)$. For each integer $1 \le i \le n$, take an East step $E:=(1,0)$ if $i \notin A$, and a North step $N:=(0,1)$ if $i \in A$.

\textit{The permutation $\pi$.} Let $x_1 < \dots < x_k$ be the elements of $A$ in increasing order. We define $\pi=\pi_1\dots \pi_k \in \mc{S}_k$ by the relation $w_i = x_{\pi_i}$ for $1 \le i \le k$.

This mapping is uniquely reversible. The path $\mb{L}$ explicitly fixes the subset $A$, the permutation $\pi$ dictates the exact arrangement of $A$ in the first $k$ positions of $w$, and the remaining $n-k$ elements are forced to monotonically fill positions $k+1$ through $n$ to satisfy $w \in \W^J$.
\end{proof}

\begin{remark}
The pairs $(\mb{L},\pi)$ can also be considered as lattice paths whose $k$ north steps are labelled sequentially by the one-line notation of the permutation $\pi$. Alternatively, these are pairs $(\la,\pi)$ where $\la$ is the partition into at most $k$ parts with largest part at most $n-k$ and whose parts are assigned distinct labels from $[k]=\set{1,\dots,k}$. The diagram of this partition can be considered as the boxes weakly above the lattice path inside the rectangle with corners at the origin and at $(n-k,k)$.
\end{remark}

\begin{example}\label{ex:min_shortest_repres_and_lattice_path}
    We illustrate the equivalent formulations given in Lemma~\ref{lemma:cosets_and_lattice_paths} with the following example. Let $w=6132457$ be a coset representative in $\mc{S}_7/\mc{S}_4$, so $n=7$ and $k=3$. First note that the $w$ is the shortest-length representative of its coset as the last $n-k =4$ entries of $w$, are increasing. The first $k=3$ positions of $w$ have the values $6,1,$ and $3$, and sorting them gives us the subset $A=\{1,3,6\}$. Thus, the corresponding lattice path from $(0,0)$ to $(4,3)$ is 
    \[ \mb{L} = NENEENE.\]
    The permutation $\pi = 312\in \mc{S}_3$ (from the ordering of $6,1,3$) labels the north steps of $\mb{L}$ sequentially, as shown in Figure~\ref{fig:lattice-path}, and the shaded region weakly above $\mb{L}$ gives the partition $\lambda = (3, 1)$. 
    \begin{figure}[htb!]
    \centering
     \begin{tikzpicture}[scale=0.85]
        \fill[darkred!25] (0,2) rectangle (3,3);
        \fill[darkred!25] (0,1) rectangle (1,2);
        \draw[step=1cm,gray,very thin] (0,0) grid (4,3);
        \draw[very thick] (0,0) -- (0,1) node[midway,left]{$3$} -- (1,1) -- (1,2) node[midway,left]{$1$} -- (3,2) -- (3,3) node[midway,left]{$2$} -- (4,3);
        \draw[fill=black] (0,0) circle (2pt) (4,3) circle (2pt);
        \node[below left] at (0,0) {$(0,0)$};
        \node[above right] at (4,3) {$(4,3)$};
    \end{tikzpicture}
    \caption{The lattice path $\mb{L}$ associated to $w=6132457$.}
    \label{fig:lattice-path}
    \end{figure}
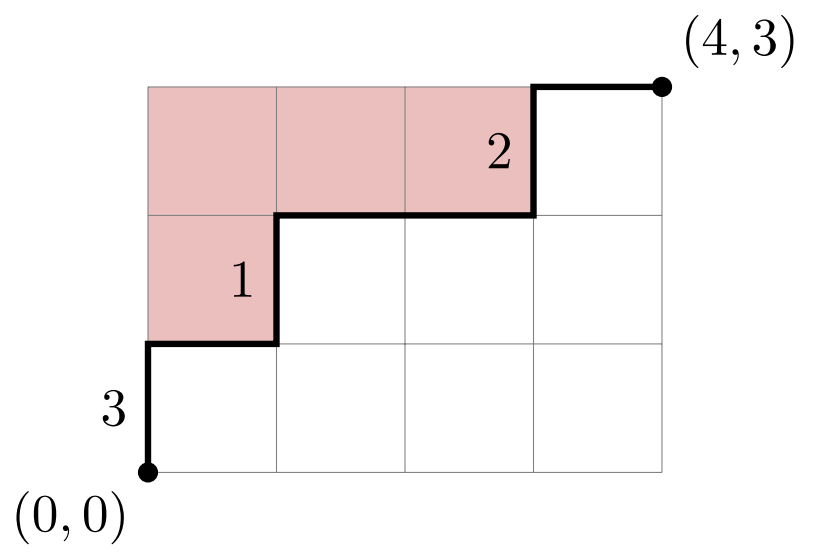
\end{example}

The proof of the next statement follows from \Cref{prop:cox_subquotients_bij} by taking $I=\set{s_1,\dots,s_{k-1}}$, $J=\set{s_{k+1},\dots,s_{n-1}}$ and $K=I\cup J$.
\begin{cor}\label{cor:order_preserving_bijection_Sn}
    There is an order-preserving bijection from the product poset $\mc{S}_n/(\mc{S}_k \x \mc{S}_{n-k}) \x \mc{S}_k$ to $(\W^J,\BC)$. That is, the multiplication map $\phi: \W^K \x \W_I \to \W^J$ given by $\phi(x,\pi) = x\pi$ is faithful.
\end{cor}

\begin{figure}[htb!]
\centering
\resizebox{\textwidth}{!}{%
\begin{tikzpicture}
\begin{scope}[shift={(0.00,0.00)},scale=1.15]
\draw[step=1,gray,very thin] (0,0) grid (3,2);
\draw[line width = 2.5pt] (0,0) -- (0,1) -- (0,2) -- (1,2) -- (2,2) -- (3,2);
\node[left] at (0.0,0.5) {$1$};
\node[left] at (0.0,1.5) {$2$};
\end{scope}
\coordinate (n0t) at (1.72,2.35);
\coordinate (n0b) at (1.72,-0.15);
\node[anchor=west] at (3.55,1.77) {$w_1=12345$};
\node[anchor=west] at (3.55,0.84) {$\pi_1=12$};

\begin{scope}[shift={(3.30,4.35)},scale=1.15]
\draw[step=1,gray,very thin] (0,0) grid (3,2);
\draw[line width = 2.5pt] (0,0) -- (0,1) -- (0,2) -- (1,2) -- (2,2) -- (3,2);
\node[left] at (0.0,0.5) {$2$};
\node[left] at (0.0,1.5) {$1$};
\end{scope}
\coordinate (n1t) at (5.02,6.70);
\coordinate (n1b) at (5.02,4.20);
\node[anchor=west] at (6.85,6.12) {$w_2=21345$};
\node[anchor=west] at (6.85,5.19) {$\pi_2=21$};

\begin{scope}[shift={(-3.30,4.35)},scale=1.15]
\fill[darkred!25] (0,1) rectangle (1,2);
\draw[step=1,gray,very thin] (0,0) grid (3,2);
\draw[line width = 2.5pt] (0,0) -- (0,1) -- (1,1) -- (1,2) -- (2,2) -- (3,2);
\node[left] at (0.0,0.5) {$1$};
\node[left] at (1.0,1.5) {$2$};
\end{scope}
\coordinate (n2t) at (-1.57,6.70);
\coordinate (n2b) at (-1.57,4.20);
\node[anchor=west] at (0.25,6.12) {$w_3=13245$};
\node[anchor=west] at (0.25,5.19) {$\pi_3=12$};

\begin{scope}[shift={(6.60,8.70)},scale=1.15]
\fill[darkred!25] (0,1) rectangle (1,2);
\draw[step=1,gray,very thin] (0,0) grid (3,2);
\draw[line width = 2.5pt] (0,0) -- (0,1) -- (1,1) -- (1,2) -- (2,2) -- (3,2);
\node[left] at (0.0,0.5) {$2$};
\node[left] at (1.0,1.5) {$1$};
\end{scope}
\coordinate (n3t) at (8.32,11.05);
\coordinate (n3b) at (8.32,8.55);
\node[anchor=west] at (10.15,10.47) {$w_4=31245$};
\node[anchor=west] at (10.15,9.54) {$\pi_4=21$};

\begin{scope}[shift={(-6.60,8.70)},scale=1.15]
\fill[darkred!25] (0,1) rectangle (2,2);
\draw[step=1,gray,very thin] (0,0) grid (3,2);
\draw[line width = 2.5pt] (0,0) -- (0,1) -- (1,1) -- (2,1) -- (2,2) -- (3,2);
\node[left] at (0.0,0.5) {$1$};
\node[left] at (2.0,1.5) {$2$};
\end{scope}
\coordinate (n4t) at (-4.88,11.05);
\coordinate (n4b) at (-4.88,8.55);
\node[anchor=west] at (-3.05,10.47) {$w_5=14235$};
\node[anchor=west] at (-3.05,9.54) {$\pi_5=12$};

\begin{scope}[shift={(0.00,8.70)},scale=1.15]
\fill[darkred!25] (0,1) rectangle (1,2);
\fill[darkred!25] (0,0) rectangle (1,1);
\draw[step=1,gray,very thin] (0,0) grid (3,2);
\draw[line width = 2.5pt] (0,0) -- (1,0) -- (1,1) -- (1,2) -- (2,2) -- (3,2);
\node[left] at (1.0,0.5) {$1$};
\node[left] at (1.0,1.5) {$2$};
\end{scope}
\coordinate (n8t) at (1.72,11.05);
\coordinate (n8b) at (1.72,8.55);
\node[anchor=west] at (3.55,10.47) {$w_6=23145$};
\node[anchor=west] at (3.55,9.54) {$\pi_6=12$};

\begin{scope}[shift={(9.90,13.05)},scale=1.15]
\fill[darkred!25] (0,1) rectangle (2,2);
\draw[step=1,gray,very thin] (0,0) grid (3,2);
\draw[line width = 2.5pt] (0,0) -- (0,1) -- (1,1) -- (2,1) -- (2,2) -- (3,2);
\node[left] at (0.0,0.5) {$2$};
\node[left] at (2.0,1.5) {$1$};
\end{scope}
\coordinate (n5t) at (11.62,15.40);
\coordinate (n5b) at (11.62,12.90);
\node[anchor=west] at (13.45,14.82) {$w_7=41235$};
\node[anchor=west] at (13.45,13.89) {$\pi_7=21$};

\begin{scope}[shift={(-9.90,13.05)},scale=1.15]
\fill[darkred!25] (0,1) rectangle (3,2);
\draw[step=1,gray,very thin] (0,0) grid (3,2);
\draw[line width = 2.5pt] (0,0) -- (0,1) -- (1,1) -- (2,1) -- (3,1) -- (3,2);
\node[left] at (0.0,0.5) {$1$};
\node[left] at (3.0,1.5) {$2$};
\end{scope}
\coordinate (n6t) at (-8.18,15.40);
\coordinate (n6b) at (-8.18,12.90);
\node[anchor=west] at (-6.35,14.82) {$w_8=15234$};
\node[anchor=west] at (-6.35,13.89) {$\pi_8=12$};

\begin{scope}[shift={(3.30,13.05)},scale=1.15]
\fill[darkred!25] (0,1) rectangle (1,2);
\fill[darkred!25] (0,0) rectangle (1,1);
\draw[step=1,gray,very thin] (0,0) grid (3,2);
\draw[line width = 2.5pt] (0,0) -- (1,0) -- (1,1) -- (1,2) -- (2,2) -- (3,2);
\node[left] at (1.0,0.5) {$2$};
\node[left] at (1.0,1.5) {$1$};
\end{scope}
\coordinate (n9t) at (5.02,15.40);
\coordinate (n9b) at (5.02,12.90);
\node[anchor=west] at (6.85,14.82) {$w_9=32145$};
\node[anchor=west] at (6.85,13.89) {$\pi_9=21$};

\begin{scope}[shift={(-3.30,13.05)},scale=1.15]
\fill[darkred!25] (0,1) rectangle (2,2);
\fill[darkred!25] (0,0) rectangle (1,1);
\draw[step=1,gray,very thin] (0,0) grid (3,2);
\draw[line width = 2.5pt] (0,0) -- (1,0) -- (1,1) -- (2,1) -- (2,2) -- (3,2);
\node[left] at (1.0,0.5) {$1$};
\node[left] at (2.0,1.5) {$2$};
\end{scope}
\coordinate (n10t) at (-1.57,15.40);
\coordinate (n10b) at (-1.57,12.90);
\node[anchor=west] at (0.25,14.82) {$w_{10}=24135$};
\node[anchor=west] at (0.25,13.89) {$\pi_{10}=12$};

\begin{scope}[shift={(9.90,17.40)},scale=1.15]
\fill[darkred!25] (0,1) rectangle (3,2);
\draw[step=1,gray,very thin] (0,0) grid (3,2);
\draw[line width = 2.5pt] (0,0) -- (0,1) -- (1,1) -- (2,1) -- (3,1) -- (3,2);
\node[left] at (0.0,0.5) {$2$};
\node[left] at (3.0,1.5) {$1$};
\end{scope}
\coordinate (n7t) at (11.62,19.75);
\coordinate (n7b) at (11.62,17.25);
\node[anchor=west] at (13.45,19.17) {$w_{11}=51234$};
\node[anchor=west] at (13.45,18.24) {$\pi_{11}=21$};

\begin{scope}[shift={(3.30,17.40)},scale=1.15]
\fill[darkred!25] (0,1) rectangle (2,2);
\fill[darkred!25] (0,0) rectangle (1,1);
\draw[step=1,gray,very thin] (0,0) grid (3,2);
\draw[line width = 2.5pt] (0,0) -- (1,0) -- (1,1) -- (2,1) -- (2,2) -- (3,2);
\node[left] at (1.0,0.5) {$2$};
\node[left] at (2.0,1.5) {$1$};
\end{scope}
\coordinate (n11t) at (5.02,19.75);
\coordinate (n11b) at (5.02,17.25);
\node[anchor=west] at (6.85,19.17) {$w_{12}=42135$};
\node[anchor=west] at (6.85,18.24) {$\pi_{12}=21$};

\begin{scope}[shift={(-9.90,17.40)},scale=1.15]
\fill[darkred!25] (0,1) rectangle (3,2);
\fill[darkred!25] (0,0) rectangle (1,1);
\draw[step=1,gray,very thin] (0,0) grid (3,2);
\draw[line width = 2.5pt] (0,0) -- (1,0) -- (1,1) -- (2,1) -- (3,1) -- (3,2);
\node[left] at (1.0,0.5) {$1$};
\node[left] at (3.0,1.5) {$2$};
\end{scope}
\coordinate (n12t) at (-8.18,19.75);
\coordinate (n12b) at (-8.18,17.25);
\node[anchor=west] at (-6.35,19.17) {$w_{13}=25134$};
\node[anchor=west] at (-6.35,18.24) {$\pi_{13}=12$};

\begin{scope}[shift={(-3.30,17.40)},scale=1.15]
\fill[darkred!25] (0,1) rectangle (2,2);
\fill[darkred!25] (0,0) rectangle (2,1);
\draw[step=1,gray,very thin] (0,0) grid (3,2);
\draw[line width = 2.5pt] (0,0) -- (1,0) -- (2,0) -- (2,1) -- (2,2) -- (3,2);
\node[left] at (2.0,0.5) {$1$};
\node[left] at (2.0,1.5) {$2$};
\end{scope}
\coordinate (n14t) at (-1.57,19.75);
\coordinate (n14b) at (-1.57,17.25);
\node[anchor=west] at (0.25,19.17) {$w_{14}=34125$};
\node[anchor=west] at (0.25,18.24) {$\pi_{14}=12$};

\begin{scope}[shift={(6.60,21.75)},scale=1.15]
\fill[darkred!25] (0,1) rectangle (3,2);
\fill[darkred!25] (0,0) rectangle (1,1);
\draw[step=1,gray,very thin] (0,0) grid (3,2);
\draw[line width = 2.5pt] (0,0) -- (1,0) -- (1,1) -- (2,1) -- (3,1) -- (3,2);
\node[left] at (1.0,0.5) {$2$};
\node[left] at (3.0,1.5) {$1$};
\end{scope}
\coordinate (n13t) at (8.32,24.10);
\coordinate (n13b) at (8.32,21.60);
\node[anchor=west] at (10.15,23.52) {$w_{15}=52134$};
\node[anchor=west] at (10.15,22.59) {$\pi_{15}=21$};

\begin{scope}[shift={(0.00,21.75)},scale=1.15]
\fill[darkred!25] (0,1) rectangle (2,2);
\fill[darkred!25] (0,0) rectangle (2,1);
\draw[step=1,gray,very thin] (0,0) grid (3,2);
\draw[line width = 2.5pt] (0,0) -- (1,0) -- (2,0) -- (2,1) -- (2,2) -- (3,2);
\node[left] at (2.0,0.5) {$2$};
\node[left] at (2.0,1.5) {$1$};
\end{scope}
\coordinate (n15t) at (1.72,24.10);
\coordinate (n15b) at (1.72,21.60);
\node[anchor=west] at (3.55,23.52) {$w_{16}=43125$};
\node[anchor=west] at (3.55,22.59) {$\pi_{16}=21$};

\begin{scope}[shift={(-6.60,21.75)},scale=1.15]
\fill[darkred!25] (0,1) rectangle (3,2);
\fill[darkred!25] (0,0) rectangle (2,1);
\draw[step=1,gray,very thin] (0,0) grid (3,2);
\draw[line width = 2.5pt] (0,0) -- (1,0) -- (2,0) -- (2,1) -- (3,1) -- (3,2);
\node[left] at (2.0,0.5) {$1$};
\node[left] at (3.0,1.5) {$2$};
\end{scope}
\coordinate (n16t) at (-4.88,24.10);
\coordinate (n16b) at (-4.88,21.60);
\node[anchor=west] at (-3.05,23.52) {$w_{17}=35124$};
\node[anchor=west] at (-3.05,22.59) {$\pi_{17}=12$};

\begin{scope}[shift={(3.30,26.10)},scale=1.15]
\fill[darkred!25] (0,1) rectangle (3,2);
\fill[darkred!25] (0,0) rectangle (2,1);
\draw[step=1,gray,very thin] (0,0) grid (3,2);
\draw[line width = 2.5pt] (0,0) -- (1,0) -- (2,0) -- (2,1) -- (3,1) -- (3,2);
\node[left] at (2.0,0.5) {$2$};
\node[left] at (3.0,1.5) {$1$};
\end{scope}
\coordinate (n17t) at (5.02,28.45);
\coordinate (n17b) at (5.02,25.95);
\node[anchor=west] at (6.85,27.87) {$w_{18}=53124$};
\node[anchor=west] at (6.85,26.94) {$\pi_{18}=21$};

\begin{scope}[shift={(-3.30,26.10)},scale=1.15]
\fill[darkred!25] (0,1) rectangle (3,2);
\fill[darkred!25] (0,0) rectangle (3,1);
\draw[step=1,gray,very thin] (0,0) grid (3,2);
\draw[line width = 2.5pt] (0,0) -- (1,0) -- (2,0) -- (3,0) -- (3,1) -- (3,2);
\node[left] at (3.0,0.5) {$1$};
\node[left] at (3.0,1.5) {$2$};
\end{scope}
\coordinate (n18t) at (-1.57,28.45);
\coordinate (n18b) at (-1.57,25.95);
\node[anchor=west] at (0.25,27.87) {$w_{19}=45123$};
\node[anchor=west] at (0.25,26.94) {$\pi_{19}=12$};

\begin{scope}[shift={(0.00,30.45)},scale=1.15]
\fill[darkred!25] (0,1) rectangle (3,2);
\fill[darkred!25] (0,0) rectangle (3,1);
\draw[step=1,gray,very thin] (0,0) grid (3,2);
\draw[line width = 2.5pt] (0,0) -- (1,0) -- (2,0) -- (3,0) -- (3,1) -- (3,2);
\node[left] at (3.0,0.5) {$2$};
\node[left] at (3.0,1.5) {$1$};
\end{scope}
\coordinate (n19t) at (1.72,32.80);
\coordinate (n19b) at (1.72,30.30);
\node[anchor=west] at (3.55,32.22) {$w_{20}=54123$};
\node[anchor=west] at (3.55,31.29) {$\pi_{20}=21$};

\draw[black] (n6t) -- (n12b);
\draw[black] (n12t) -- (n16b);
\draw[black] (n4t) -- (n6b);
\draw[black] (n12t) -- (n13b);
\draw[black] (n5t) -- (n7b);
\draw[black] (n14t) -- (n16b);
\draw[black] (n0t) -- (n2b);
\draw[black] (n8t) -- (n9b);
\draw[black] (n9t) -- (n11b);
\draw[black] (n1t) -- (n3b);
\draw[black] (n10t) -- (n12b);
\draw[black] (n2t) -- (n8b);
\draw[black] (n13t) -- (n17b);
\draw[black] (n15t) -- (n17b);
\draw[black] (n7t) -- (n13b);
\draw[black] (n18t) -- (n19b);
\draw[black] (n4t) -- (n5b);
\draw[black] (n3t) -- (n9b);
\draw[black] (n14t) -- (n15b);
\draw[black] (n0t) -- (n1b);
\draw[black] (n2t) -- (n4b);
\draw[black] (n10t) -- (n11b);
\draw[black] (n10t) -- (n14b);
\draw[black] (n11t) -- (n13b);
\draw[black] (n6t) -- (n7b);
\draw[black] (n16t) -- (n18b);
\draw[black] (n3t) -- (n5b);
\draw[black] (n4t) -- (n10b);
\draw[black] (n5t) -- (n11b);
\draw[black] (n8t) -- (n10b);
\draw[black] (n17t) -- (n19b);
\draw[black] (n2t) -- (n3b);
\draw[black] (n11t) -- (n15b);
\draw[black] (n16t) -- (n17b);
\draw[red,line width=3pt] (n1t) -- (n8b);
\draw[red,line width=3pt] (n15t) -- (n18b);
\draw[red,line width=3pt] (n9t) -- (n14b);
\end{tikzpicture}
}
\caption{The bijection $\phi:L(2,3)\x \mc{S}_2 \to \mc{S}_5/\mc{S}_3$. Black edges are covering relations shared by both orders; red edges are covering relations present in the Bruhat--Chevalley order but not in the product order.}
\label{fig:lattice-hasse-full}
\end{figure}
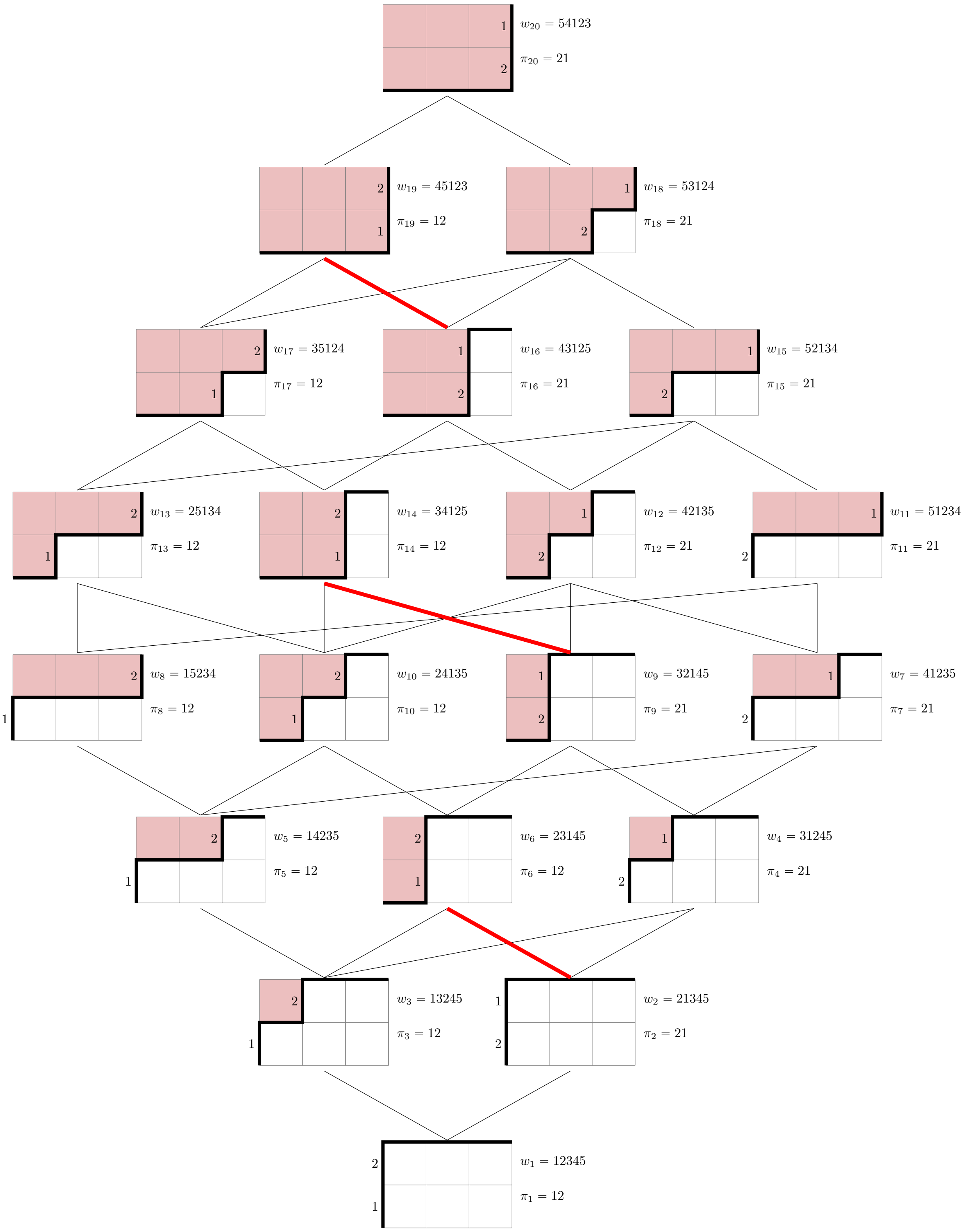

\begin{example}\label{ex:application_preserving_bijection_Sn}
Let $\W = \mc{S}_5$ and $k=2$, so that $J = \{s_3, s_4\}$ and $\W^J \cong \mc{S}_5/\mc{S}_3$. Take 
$I = \{s_1\}$ and $K = I \cup J= \{s_1, s_3, s_4\}$. Then $\W^K \cong \mc{S}_5/(\mc{S}_2 \times \mc{S}_3)$ and $\W_I \cong \mc{S}_2$. Consider the pair $(y, \tau) \in \W^K \x \W_I$ with
\begin{itemize}
    \item $y = 34125$. Since it has reduced word $s_2 s_1 s_3 s_2$, then $\ell(y) = 4$.
    \item $\tau = 21345$. This acts only on the first two positions. Its reduced word is $s_1$, so $\ell(\tau) = 1$.
\end{itemize}
Applying the bijection yields $w = \phi(y, \tau) = y\tau = 43125$. 
Because $\phi$ is length-additive, the reduced word for $w$ is the direct concatenation:
\[ w = (s_2 s_1 s_3 s_2)(s_1) = s_2 s_1 s_3 s_2 s_1. \]
\noindent
To see that $\phi$ is order-preserving, we introduce another pair $(x, \pi) \in \W^K \x \W_I$ where:
\begin{itemize}
    \item $x = 23145$, with reduced word $s_1 s_2$.  Note that $x \leq_{BC} y$ because $s_1 s_2$ is a subexpression of $s_2 s_1 s_3 s_2$.
    \item $\pi = 12345$ (the identity $e$). Note that $\pi \leq_{BC} \tau$.
\end{itemize}
Since inequality holds in both coordinates independently, $(x, \pi) \leq (y, \tau)$ in the Cartesian product order.

We now apply the bijection to this pair, yielding $z = \phi(x, \pi) = (23145)(e) = 23145$. 
Its reduced word is simply $s_1 s_2$.  Comparing the two mapped elements in the quotient $\W^J$:
\begin{align*}
    w &= s_2 s_1 s_3 s_2 s_1, \\
    z &= s_1 s_2. 
\end{align*}
It is immediately clear that the reduced word for $z$ is a subexpression of the reduced word for $w$. By \Cref{thm:bruhat-subexp}, this explicitly proves that $z \leq_{BC} w$, verifying that the product order embeds faithfully into the quotient order as illustrated in Figure~\ref{fig:lattice-hasse-full}.
    \end{example}

\subsection{Codes and SCDs} \label{sec:results.A}
Continuing with the notation of the previous subsection, here we show that $(\W^J,\BC)$ is an SCO. The proof comes by restricting the code for $\mc{S}_n$ to the set $\W^J$ of shortest length coset representatives for $\mc{S}_n/\mc{S}_{n-k}$.  We begin by recalling the $\mc{S}_n$ statement, attributed as an observation due to E. Gansner in \cite[p.~182]{stanley1980}.
\begin{prop}  \label{prop:Sn_SCD}
    The bijection $\code^{-1}:\mb{n}_n \to \mc{S}_n$ is an order-preserving poset map between the product order on chains and the Bruhat--Chevalley order on $\mc{S}_n$. 
\end{prop}
\begin{proof}
    The bijection follows from \Cref{lem:staircase bijection}. Assume that $\alpha\preceq \beta$ for some $\alpha,\beta\in \mb{n}_n$ and consider $w(\alpha)$ and $w(\beta)$ with $\ul{{w}(\alpha)}=\ul{r}_\alpha(1)\cdots\ul{r}_\alpha{(n-1)}$ and $\ul{{w}(\beta)}=\ul{r}_\beta(1)\cdots \ul{r}_\beta{(n-1)}$. Observe that the inequality of compositions $\alpha\preceq \beta$ in $\mb{n}_n$ implies that $a_i\leq b_i$ for each $i$. Thus, $\ul{r}_\alpha(i)$ is a subexpression of $\ul{r}_\beta(i)$ for each $i$. By the subexpression description of Bruhat--Chevalley order given in Theorem~\ref{thm:bruhat-subexp}, it is then clear that $w(\alpha)\leq_{BC} w(\beta)$.
\end{proof}
The following corollary is then a direct consequence of \Cref{prop:Sn_SCD} and  \Cref{lem:product_SCD}.
\begin{cor}[{\cite[p. 182]{stanley1980}}]
    The poset $(\mc{S}_n, \leq_{BC})$ is an SCO, with SCDs provided by passing any SCD of $\mb{n}_n$ through the map $\code^{-1}$.
\end{cor}
Our main result for this subsection follows by showing that the usual code on $\mc{S}_n$ restricts nicely to the quotients $\W^J$, providing SCDs. 

\begin{lemma} \label{lem:restricted_code}
    For any $w\in \W^J$, $\code(w)\in \mb{n}_k$. 
\end{lemma}
\begin{proof}
    We argue that $\code(w)=(a_1,\dots,a_n)\in \mb{n}_k$ for any $w \in \W^J$.  By Lemma~\ref{lem:min_reps}, we have $w_{k+1}<\dots<w_n$ in the one-line notation for $w$. This implies that for columns $k+1,\dots, n$ there are no boxes in row $i$ of the ray diagram of $w$ for any $i\in \set{w_{k+1},\dots, w_n}$. But since $\set{w_1,\dots,w_k}=[n]-\set{w_{k+1},\dots, w_{k}}$, the ray diagram construction implies that there are no boxes in these columns for any row at all. Thus, $a_{k+1}=a_{k+2}=\dots=a_n=0$, and the composition $\code(w)$ fits inside the truncated staircase shape determined by $\mb{n}_k$.
\end{proof}
\begin{theorem}\label{thm:WJ_SCO}
    There is an order-preserving bijection $\code^{-1}:\mb{n}_k \to \W^J$. Consequently, $(\W^J,\BC)$ is an SCO.
\end{theorem}
\begin{proof}
It is clear that both sets have size $(n)_k$. Since, by \Cref{lem:restricted_code}, the restriction of the code map to $\W^J$ is an injection into an equinumerous set, it is a bijection. Hence, so is its inverse $\code^{-1}:\mb{n}_k \to \W^J$. 

Supposing that $\code(w) \preceq \code(z)$ for $w,z\in \W^J$, \Cref{prop:Sn_SCD} implies that $w \leq_{BC} z$ as elements of $\mc{S}_n$. Then the fact that the projection map $P^J: \W\to \W^J$ is order-preserving \cite[Proposition 2.5.1]{bjornerBrenti} implies that $\code^{-1}:\mb{n}_k \to \W^J$ is also order-preserving.
\end{proof}  
\begin{remark}
    Observe that if $L(k,n-k)$ were known to be an SCO, then \Cref{cor:order_preserving_bijection_Sn}, \Cref{prop:Sn_SCD}, and \Cref{lem:product_SCD} would combine to provide an alternative proof of the SCO statement in \Cref{thm:WJ_SCO}. 
\end{remark}
\begin{remark} \label{rk:Sn-k}
    \Cref{thm:WJ_SCO} generalizes to give a bijection $\mb{n}_k\to \W^J$ for any $J\sus S$ such that $\W_J$ is irreducible. That is, one can also take $\W_J \cong \mc{S}_{n-k}$ where $J=\set{s_i,\dots,s_{i+n-k-2}}$ for some other $i<k+1$. However, the order-preserving bijection is no longer given by the restriction of $\code^{-1}:\mb{n}_n\to\mc{S}_n$ as described here, but rather a composite construction. This is is due to the need to take both left and right parabolic quotients, under the perspective of \Cref{rk:code.BPseq}. We omit details as this form of quotient is less directly relevant to $L(k,n-k)$.
\end{remark}
\begin{example}
We illustrate an instance of symmetric chain decompositions induced by the poset map $\code^{-1}:\ul{5}\x \ul{4}\to \mc{S}_{5}/\mc{S}_3$ in Figure~\ref{fig:SCD}. 
\begin{figure}[htbp]
\centering
\begin{subfigure}[b]{0.49\textwidth}
\centering
\resizebox{\linewidth}{!}{%
\begin{tikzpicture}[scale=1]
\begin{scope}[shift={(0.00,0.00)},scale=0.75]
\draw[step=1,gray,very thin] (0,0) grid (3,2);
\draw[very thick] (0,0) -- (0,1) -- (0,2) -- (1,2) -- (2,2) -- (3,2);
\end{scope}
\coordinate (n0t) at (1.05,1.55);
\coordinate (n0b) at (1.05,-0.15);

\begin{scope}[shift={(3.30,4.35)},scale=0.75]
\draw[step=1,gray,very thin] (0,0) grid (3,2);
\draw[very thick] (0,0) -- (0,1) -- (0,2) -- (1,2) -- (2,2) -- (3,2);
\end{scope}
\coordinate (n1t) at (4.35,5.90);
\coordinate (n1b) at (4.35,4.20);

\begin{scope}[shift={(-3.30,4.35)},scale=0.75]
\fill[darkred!25] (0,1) rectangle (1,2);
\draw[step=1,gray,very thin] (0,0) grid (3,2);
\draw[very thick] (0,0) -- (0,1) -- (1,1) -- (1,2) -- (2,2) -- (3,2);
\end{scope}
\coordinate (n2t) at (-2.25,5.90);
\coordinate (n2b) at (-2.25,4.20);

\begin{scope}[shift={(6.60,8.70)},scale=0.75]
\fill[darkred!25] (0,1) rectangle (1,2);
\draw[step=1,gray,very thin] (0,0) grid (3,2);
\draw[very thick] (0,0) -- (0,1) -- (1,1) -- (1,2) -- (2,2) -- (3,2);
\end{scope}
\coordinate (n3t) at (7.65,10.25);
\coordinate (n3b) at (7.65,8.55);

\begin{scope}[shift={(-6.60,8.70)},scale=0.75]
\fill[darkred!25] (0,1) rectangle (2,2);
\draw[step=1,gray,very thin] (0,0) grid (3,2);
\draw[very thick] (0,0) -- (0,1) -- (1,1) -- (2,1) -- (2,2) -- (3,2);
\end{scope}
\coordinate (n4t) at (-5.55,10.25);
\coordinate (n4b) at (-5.55,8.55);

\begin{scope}[shift={(0.00,8.70)},scale=0.75]
\fill[darkred!25] (0,1) rectangle (1,2);
\fill[darkred!25] (0,0) rectangle (1,1);
\draw[step=1,gray,very thin] (0,0) grid (3,2);
\draw[very thick] (0,0) -- (1,0) -- (1,1) -- (1,2) -- (2,2) -- (3,2);
\end{scope}
\coordinate (n8t) at (1.05,10.25);
\coordinate (n8b) at (1.05,8.55);

\begin{scope}[shift={(9.90,13.05)},scale=0.75]
\fill[darkred!25] (0,1) rectangle (2,2);
\draw[step=1,gray,very thin] (0,0) grid (3,2);
\draw[very thick] (0,0) -- (0,1) -- (1,1) -- (2,1) -- (2,2) -- (3,2);
\end{scope}
\coordinate (n5t) at (10.95,14.60);
\coordinate (n5b) at (10.95,12.90);

\begin{scope}[shift={(-9.90,13.05)},scale=0.75]
\fill[darkred!25] (0,1) rectangle (3,2);
\draw[step=1,gray,very thin] (0,0) grid (3,2);
\draw[very thick] (0,0) -- (0,1) -- (1,1) -- (2,1) -- (3,1) -- (3,2);
\end{scope}
\coordinate (n6t) at (-8.85,14.60);
\coordinate (n6b) at (-8.85,12.90);

\begin{scope}[shift={(3.30,13.05)},scale=0.75]
\fill[darkred!25] (0,1) rectangle (1,2);
\fill[darkred!25] (0,0) rectangle (1,1);
\draw[step=1,gray,very thin] (0,0) grid (3,2);
\draw[very thick] (0,0) -- (1,0) -- (1,1) -- (1,2) -- (2,2) -- (3,2);
\end{scope}
\coordinate (n9t) at (4.35,14.60);
\coordinate (n9b) at (4.35,12.90);

\begin{scope}[shift={(-3.30,13.05)},scale=0.75]
\fill[darkred!25] (0,1) rectangle (2,2);
\fill[darkred!25] (0,0) rectangle (1,1);
\draw[step=1,gray,very thin] (0,0) grid (3,2);
\draw[very thick] (0,0) -- (1,0) -- (1,1) -- (2,1) -- (2,2) -- (3,2);
\end{scope}
\coordinate (n10t) at (-2.25,14.60);
\coordinate (n10b) at (-2.25,12.90);

\begin{scope}[shift={(9.90,17.40)},scale=0.75]
\fill[darkred!25] (0,1) rectangle (3,2);
\draw[step=1,gray,very thin] (0,0) grid (3,2);
\draw[very thick] (0,0) -- (0,1) -- (1,1) -- (2,1) -- (3,1) -- (3,2);
\end{scope}
\coordinate (n7t) at (10.95,18.95);
\coordinate (n7b) at (10.95,17.25);

\begin{scope}[shift={(3.30,17.40)},scale=0.75]
\fill[darkred!25] (0,1) rectangle (2,2);
\fill[darkred!25] (0,0) rectangle (1,1);
\draw[step=1,gray,very thin] (0,0) grid (3,2);
\draw[very thick] (0,0) -- (1,0) -- (1,1) -- (2,1) -- (2,2) -- (3,2);
\end{scope}
\coordinate (n11t) at (4.35,18.95);
\coordinate (n11b) at (4.35,17.25);

\begin{scope}[shift={(-9.90,17.40)},scale=0.75]
\fill[darkred!25] (0,1) rectangle (3,2);
\fill[darkred!25] (0,0) rectangle (1,1);
\draw[step=1,gray,very thin] (0,0) grid (3,2);
\draw[very thick] (0,0) -- (1,0) -- (1,1) -- (2,1) -- (3,1) -- (3,2);
\end{scope}
\coordinate (n12t) at (-8.85,18.95);
\coordinate (n12b) at (-8.85,17.25);

\begin{scope}[shift={(-3.30,17.40)},scale=0.75]
\fill[darkred!25] (0,1) rectangle (2,2);
\fill[darkred!25] (0,0) rectangle (2,1);
\draw[step=1,gray,very thin] (0,0) grid (3,2);
\draw[very thick] (0,0) -- (1,0) -- (2,0) -- (2,1) -- (2,2) -- (3,2);
\end{scope}
\coordinate (n14t) at (-2.25,18.95);
\coordinate (n14b) at (-2.25,17.25);

\begin{scope}[shift={(6.60,21.75)},scale=0.75]
\fill[darkred!25] (0,1) rectangle (3,2);
\fill[darkred!25] (0,0) rectangle (1,1);
\draw[step=1,gray,very thin] (0,0) grid (3,2);
\draw[very thick] (0,0) -- (1,0) -- (1,1) -- (2,1) -- (3,1) -- (3,2);
\end{scope}
\coordinate (n13t) at (7.65,23.30);
\coordinate (n13b) at (7.65,21.60);

\begin{scope}[shift={(0.00,21.75)},scale=0.75]
\fill[darkred!25] (0,1) rectangle (2,2);
\fill[darkred!25] (0,0) rectangle (2,1);
\draw[step=1,gray,very thin] (0,0) grid (3,2);
\draw[very thick] (0,0) -- (1,0) -- (2,0) -- (2,1) -- (2,2) -- (3,2);
\end{scope}
\coordinate (n15t) at (1.05,23.30);
\coordinate (n15b) at (1.05,21.60);

\begin{scope}[shift={(-6.60,21.75)},scale=0.75]
\fill[darkred!25] (0,1) rectangle (3,2);
\fill[darkred!25] (0,0) rectangle (2,1);
\draw[step=1,gray,very thin] (0,0) grid (3,2);
\draw[very thick] (0,0) -- (1,0) -- (2,0) -- (2,1) -- (3,1) -- (3,2);
\end{scope}
\coordinate (n16t) at (-5.55,23.30);
\coordinate (n16b) at (-5.55,21.60);

\begin{scope}[shift={(3.30,26.10)},scale=0.75]
\fill[darkred!25] (0,1) rectangle (3,2);
\fill[darkred!25] (0,0) rectangle (2,1);
\draw[step=1,gray,very thin] (0,0) grid (3,2);
\draw[very thick] (0,0) -- (1,0) -- (2,0) -- (2,1) -- (3,1) -- (3,2);
\end{scope}
\coordinate (n17t) at (4.35,27.65);
\coordinate (n17b) at (4.35,25.95);

\begin{scope}[shift={(-3.30,26.10)},scale=0.75]
\fill[darkred!25] (0,1) rectangle (3,2);
\fill[darkred!25] (0,0) rectangle (3,1);
\draw[step=1,gray,very thin] (0,0) grid (3,2);
\draw[very thick] (0,0) -- (1,0) -- (2,0) -- (3,0) -- (3,1) -- (3,2);
\end{scope}
\coordinate (n18t) at (-2.25,27.65);
\coordinate (n18b) at (-2.25,25.95);

\begin{scope}[shift={(0.00,30.45)},scale=0.75]
\fill[darkred!25] (0,1) rectangle (3,2);
\fill[darkred!25] (0,0) rectangle (3,1);
\draw[step=1,gray,very thin] (0,0) grid (3,2);
\draw[very thick] (0,0) -- (1,0) -- (2,0) -- (3,0) -- (3,1) -- (3,2);
\end{scope}
\coordinate (n19t) at (1.05,32.00);
\coordinate (n19b) at (1.05,30.30);

\draw[gray!40] (n10t) -- (n11b);
\draw[gray!40] (n10t) -- (n12b);
\draw[gray!40] (n5t) -- (n11b);
\draw[gray!40] (n12t) -- (n13b);
\draw[gray!40] (n9t) -- (n14b);
\draw[gray!40] (n16t) -- (n17b);
\draw[gray!40] (n2t) -- (n3b);
\draw[gray!40] (n4t) -- (n5b);
\draw[gray!40] (n6t) -- (n7b);
\draw[gray!40] (n0t) -- (n2b);
\draw[gray!40] (n1t) -- (n8b);
\draw[gray!40] (n8t) -- (n10b);
\draw[gray!40] (n11t) -- (n13b);
\draw[gray!40] (n12t) -- (n16b);
\draw[gray!40] (n14t) -- (n15b);
\draw[gray!40] (n15t) -- (n17b);
\draw[gray!40] (n16t) -- (n18b);
\draw[gray!40] (n18t) -- (n19b);
\draw[gray!40] (n2t) -- (n4b);
\draw[gray!40] (n3t) -- (n9b);
\draw[gray!40] (n4t) -- (n6b);

\draw[blue,line width=6pt] (n0t) -- (n1b);
\draw[blue,line width=6pt] (n1t) -- (n3b);
\draw[blue,line width=6pt] (n3t) -- (n5b);
\draw[blue,line width=6pt] (n5t) -- (n7b);
\draw[blue,line width=6pt] (n7t) -- (n13b);
\draw[blue,line width=6pt] (n13t) -- (n17b);
\draw[blue,line width=6pt] (n17t) -- (n19b);

\draw[ForestGreen,line width=6pt] (n2t) -- (n8b);
\draw[ForestGreen,line width=6pt] (n8t) -- (n9b);
\draw[ForestGreen,line width=6pt] (n9t) -- (n11b);
\draw[ForestGreen,line width=6pt] (n11t) -- (n15b);
\draw[ForestGreen,line width=6pt] (n15t) -- (n18b);

\draw[orange,line width=6pt] (n4t) -- (n10b);
\draw[orange,line width=6pt] (n10t) -- (n14b);
\draw[orange,line width=6pt] (n14t) -- (n16b);

\draw[violet,line width=6pt] (n6t) -- (n12b);

\end{tikzpicture}%
}
\caption{}
\label{fig:SCD_a}
\end{subfigure}%
\hfill
\begin{subfigure}[b]{0.49\textwidth}
\centering
\resizebox{\linewidth}{!}{%
\begin{tikzpicture}[scale=1]
\begin{scope}[shift={(0.00,0.00)},scale=0.75]
\draw[step=1,gray,very thin] (0,0) grid (3,2);
\draw[very thick] (0,0) -- (0,1) -- (0,2) -- (1,2) -- (2,2) -- (3,2);
\end{scope}
\coordinate (n0t) at (1.05,1.55);
\coordinate (n0b) at (1.05,-0.15);

\begin{scope}[shift={(3.30,4.35)},scale=0.75]
\draw[step=1,gray,very thin] (0,0) grid (3,2);
\draw[very thick] (0,0) -- (0,1) -- (0,2) -- (1,2) -- (2,2) -- (3,2);
\end{scope}
\coordinate (n1t) at (4.35,5.90);
\coordinate (n1b) at (4.35,4.20);

\begin{scope}[shift={(-3.30,4.35)},scale=0.75]
\fill[darkred!25] (0,1) rectangle (1,2);
\draw[step=1,gray,very thin] (0,0) grid (3,2);
\draw[very thick] (0,0) -- (0,1) -- (1,1) -- (1,2) -- (2,2) -- (3,2);
\end{scope}
\coordinate (n2t) at (-2.25,5.90);
\coordinate (n2b) at (-2.25,4.20);

\begin{scope}[shift={(6.60,8.70)},scale=0.75]
\fill[darkred!25] (0,1) rectangle (1,2);
\draw[step=1,gray,very thin] (0,0) grid (3,2);
\draw[very thick] (0,0) -- (0,1) -- (1,1) -- (1,2) -- (2,2) -- (3,2);
\end{scope}
\coordinate (n3t) at (7.65,10.25);
\coordinate (n3b) at (7.65,8.55);

\begin{scope}[shift={(-6.60,8.70)},scale=0.75]
\fill[darkred!25] (0,1) rectangle (2,2);
\draw[step=1,gray,very thin] (0,0) grid (3,2);
\draw[very thick] (0,0) -- (0,1) -- (1,1) -- (2,1) -- (2,2) -- (3,2);
\end{scope}
\coordinate (n4t) at (-5.55,10.25);
\coordinate (n4b) at (-5.55,8.55);

\begin{scope}[shift={(0.00,8.70)},scale=0.75]
\fill[darkred!25] (0,1) rectangle (1,2);
\fill[darkred!25] (0,0) rectangle (1,1);
\draw[step=1,gray,very thin] (0,0) grid (3,2);
\draw[very thick] (0,0) -- (1,0) -- (1,1) -- (1,2) -- (2,2) -- (3,2);
\end{scope}
\coordinate (n8t) at (1.05,10.25);
\coordinate (n8b) at (1.05,8.55);

\begin{scope}[shift={(9.90,13.05)},scale=0.75]
\fill[darkred!25] (0,1) rectangle (2,2);
\draw[step=1,gray,very thin] (0,0) grid (3,2);
\draw[very thick] (0,0) -- (0,1) -- (1,1) -- (2,1) -- (2,2) -- (3,2);
\end{scope}
\coordinate (n5t) at (10.95,14.60);
\coordinate (n5b) at (10.95,12.90);

\begin{scope}[shift={(-9.90,13.05)},scale=0.75]
\fill[darkred!25] (0,1) rectangle (3,2);
\draw[step=1,gray,very thin] (0,0) grid (3,2);
\draw[very thick] (0,0) -- (0,1) -- (1,1) -- (2,1) -- (3,1) -- (3,2);
\end{scope}
\coordinate (n6t) at (-8.85,14.60);
\coordinate (n6b) at (-8.85,12.90);

\begin{scope}[shift={(3.30,13.05)},scale=0.75]
\fill[darkred!25] (0,1) rectangle (1,2);
\fill[darkred!25] (0,0) rectangle (1,1);
\draw[step=1,gray,very thin] (0,0) grid (3,2);
\draw[very thick] (0,0) -- (1,0) -- (1,1) -- (1,2) -- (2,2) -- (3,2);
\end{scope}
\coordinate (n9t) at (4.35,14.60);
\coordinate (n9b) at (4.35,12.90);

\begin{scope}[shift={(-3.30,13.05)},scale=0.75]
\fill[darkred!25] (0,1) rectangle (2,2);
\fill[darkred!25] (0,0) rectangle (1,1);
\draw[step=1,gray,very thin] (0,0) grid (3,2);
\draw[very thick] (0,0) -- (1,0) -- (1,1) -- (2,1) -- (2,2) -- (3,2);
\end{scope}
\coordinate (n10t) at (-2.25,14.60);
\coordinate (n10b) at (-2.25,12.90);

\begin{scope}[shift={(9.90,17.40)},scale=0.75]
\fill[darkred!25] (0,1) rectangle (3,2);
\draw[step=1,gray,very thin] (0,0) grid (3,2);
\draw[very thick] (0,0) -- (0,1) -- (1,1) -- (2,1) -- (3,1) -- (3,2);
\end{scope}
\coordinate (n7t) at (10.95,18.95);
\coordinate (n7b) at (10.95,17.25);

\begin{scope}[shift={(3.30,17.40)},scale=0.75]
\fill[darkred!25] (0,1) rectangle (2,2);
\fill[darkred!25] (0,0) rectangle (1,1);
\draw[step=1,gray,very thin] (0,0) grid (3,2);
\draw[very thick] (0,0) -- (1,0) -- (1,1) -- (2,1) -- (2,2) -- (3,2);
\end{scope}
\coordinate (n11t) at (4.35,18.95);
\coordinate (n11b) at (4.35,17.25);

\begin{scope}[shift={(-9.90,17.40)},scale=0.75]
\fill[darkred!25] (0,1) rectangle (3,2);
\fill[darkred!25] (0,0) rectangle (1,1);
\draw[step=1,gray,very thin] (0,0) grid (3,2);
\draw[very thick] (0,0) -- (1,0) -- (1,1) -- (2,1) -- (3,1) -- (3,2);
\end{scope}
\coordinate (n12t) at (-8.85,18.95);
\coordinate (n12b) at (-8.85,17.25);

\begin{scope}[shift={(-3.30,17.40)},scale=0.75]
\fill[darkred!25] (0,1) rectangle (2,2);
\fill[darkred!25] (0,0) rectangle (2,1);
\draw[step=1,gray,very thin] (0,0) grid (3,2);
\draw[very thick] (0,0) -- (1,0) -- (2,0) -- (2,1) -- (2,2) -- (3,2);
\end{scope}
\coordinate (n14t) at (-2.25,18.95);
\coordinate (n14b) at (-2.25,17.25);

\begin{scope}[shift={(6.60,21.75)},scale=0.75]
\fill[darkred!25] (0,1) rectangle (3,2);
\fill[darkred!25] (0,0) rectangle (1,1);
\draw[step=1,gray,very thin] (0,0) grid (3,2);
\draw[very thick] (0,0) -- (1,0) -- (1,1) -- (2,1) -- (3,1) -- (3,2);
\end{scope}
\coordinate (n13t) at (7.65,23.30);
\coordinate (n13b) at (7.65,21.60);

\begin{scope}[shift={(0.00,21.75)},scale=0.75]
\fill[darkred!25] (0,1) rectangle (2,2);
\fill[darkred!25] (0,0) rectangle (2,1);
\draw[step=1,gray,very thin] (0,0) grid (3,2);
\draw[very thick] (0,0) -- (1,0) -- (2,0) -- (2,1) -- (2,2) -- (3,2);
\end{scope}
\coordinate (n15t) at (1.05,23.30);
\coordinate (n15b) at (1.05,21.60);

\begin{scope}[shift={(-6.60,21.75)},scale=0.75]
\fill[darkred!25] (0,1) rectangle (3,2);
\fill[darkred!25] (0,0) rectangle (2,1);
\draw[step=1,gray,very thin] (0,0) grid (3,2);
\draw[very thick] (0,0) -- (1,0) -- (2,0) -- (2,1) -- (3,1) -- (3,2);
\end{scope}
\coordinate (n16t) at (-5.55,23.30);
\coordinate (n16b) at (-5.55,21.60);

\begin{scope}[shift={(3.30,26.10)},scale=0.75]
\fill[darkred!25] (0,1) rectangle (3,2);
\fill[darkred!25] (0,0) rectangle (2,1);
\draw[step=1,gray,very thin] (0,0) grid (3,2);
\draw[very thick] (0,0) -- (1,0) -- (2,0) -- (2,1) -- (3,1) -- (3,2);
\end{scope}
\coordinate (n17t) at (4.35,27.65);
\coordinate (n17b) at (4.35,25.95);

\begin{scope}[shift={(-3.30,26.10)},scale=0.75]
\fill[darkred!25] (0,1) rectangle (3,2);
\fill[darkred!25] (0,0) rectangle (3,1);
\draw[step=1,gray,very thin] (0,0) grid (3,2);
\draw[very thick] (0,0) -- (1,0) -- (2,0) -- (3,0) -- (3,1) -- (3,2);
\end{scope}
\coordinate (n18t) at (-2.25,27.65);
\coordinate (n18b) at (-2.25,25.95);

\begin{scope}[shift={(0.00,30.45)},scale=0.75]
\fill[darkred!25] (0,1) rectangle (3,2);
\fill[darkred!25] (0,0) rectangle (3,1);
\draw[step=1,gray,very thin] (0,0) grid (3,2);
\draw[very thick] (0,0) -- (1,0) -- (2,0) -- (3,0) -- (3,1) -- (3,2);
\end{scope}
\coordinate (n19t) at (1.05,32.00);
\coordinate (n19b) at (1.05,30.30);

\draw[gray!40] (n10t) -- (n11b);
\draw[gray!40] (n10t) -- (n12b);
\draw[gray!40] (n5t) -- (n11b);
\draw[gray!40] (n12t) -- (n13b);
\draw[gray!40] (n9t) -- (n14b);
\draw[gray!40] (n16t) -- (n17b);
\draw[gray!40] (n2t) -- (n3b);
\draw[gray!40] (n4t) -- (n5b);
\draw[gray!40] (n6t) -- (n7b);
\draw[gray!40] (n0t) -- (n1b);
\draw[gray!40] (n1t) -- (n3b);
\draw[gray!40] (n3t) -- (n5b);
\draw[gray!40] (n7t) -- (n13b);
\draw[gray!40] (n13t) -- (n17b);
\draw[gray!40] (n17t) -- (n19b);
\draw[gray!40] (n2t) -- (n8b);
\draw[gray!40] (n8t) -- (n9b);
\draw[gray!40] (n11t) -- (n15b);
\draw[gray!40] (n15t) -- (n18b);
\draw[gray!40] (n4t) -- (n10b);
\draw[gray!40] (n14t) -- (n16b);

\draw[blue,line width=6pt] (n0t) -- (n2b);
\draw[blue,line width=6pt] (n2t) -- (n4b);
\draw[blue,line width=6pt] (n4t) -- (n6b);
\draw[blue,line width=6pt] (n6t) -- (n12b);
\draw[blue,line width=6pt] (n12t) -- (n16b);
\draw[blue,line width=6pt] (n16t) -- (n18b);
\draw[blue,line width=6pt] (n18t) -- (n19b);

\draw[ForestGreen,line width=6pt] (n1t) -- (n8b);
\draw[ForestGreen,line width=6pt] (n8t) -- (n10b);
\draw[ForestGreen,line width=6pt] (n10t) -- (n14b);
\draw[ForestGreen,line width=6pt] (n14t) -- (n15b);
\draw[ForestGreen,line width=6pt] (n15t) -- (n17b);

\draw[orange,line width=6pt] (n3t) -- (n9b);
\draw[orange,line width=6pt] (n9t) -- (n11b);
\draw[orange,line width=6pt] (n11t) -- (n13b);

\draw[violet,line width=6pt] (n5t) -- (n7b);

\end{tikzpicture}%
}
\caption{}
\label{fig:SCD_b}
\end{subfigure}
\caption{Two symmetric chain decompositions: (A) an SCD which ascend in the first coordinate of $\ul{5}\x \ul{4}$ first and (B) an SCD which ascends in the second coordinate first. Observe that each SCD uses some edge which is not present in the product order $L(2,3) \x \mc{S}_2$ of \Cref{lem:restricted_code} (see \Cref{fig:lattice-hasse-full}).}
\label{fig:SCD}
\end{figure}
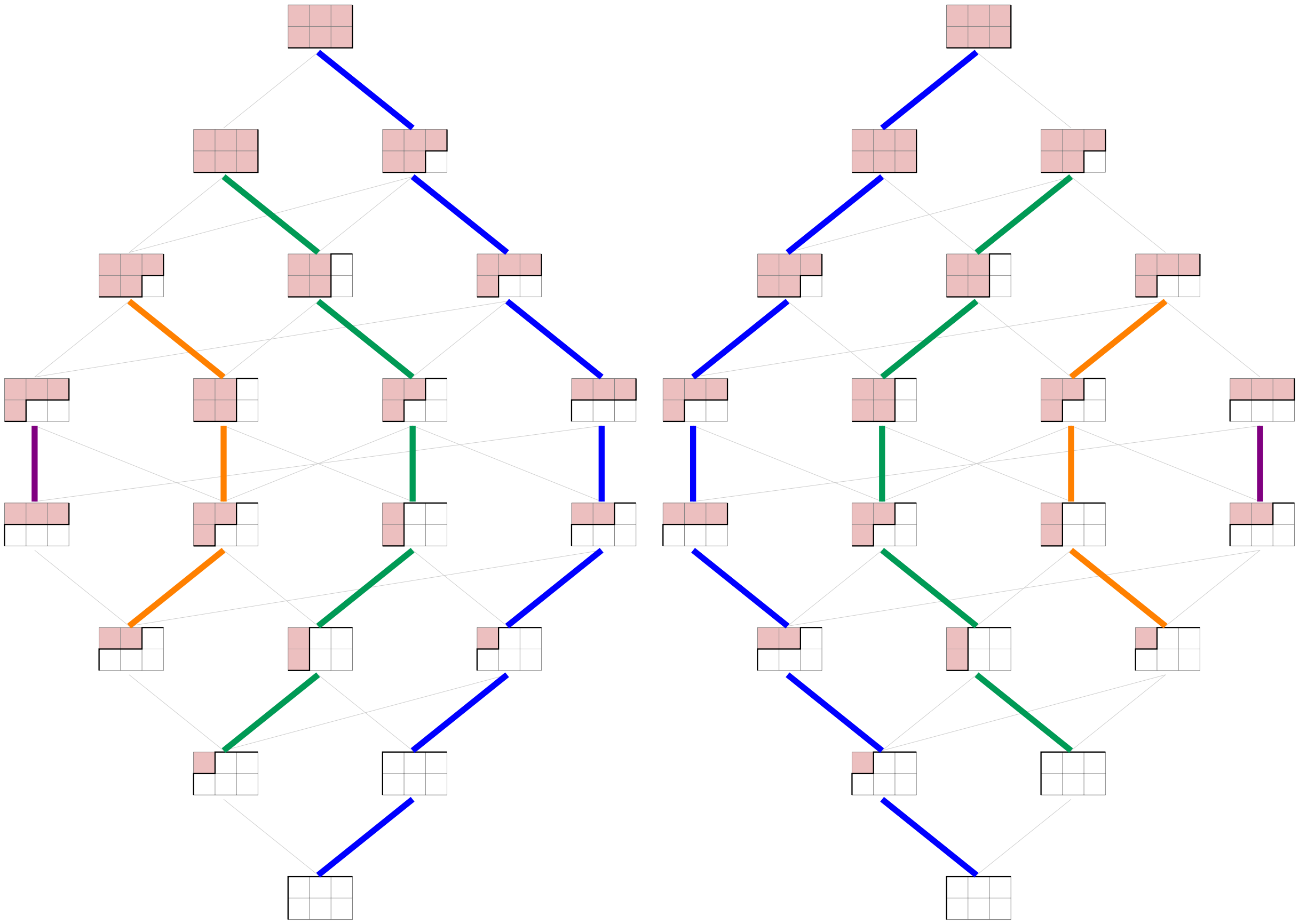
\end{example}

\subsection{Weak compositions inside a rectangle}
Analogous to partitions, (weak) compositions can be ordered by inclusion of diagrams to create a poset. As with Young's lattice, restriction of this poset an initial interval $[\varnothing, \alpha]$ leads to numerous interesting combinatorial questions. For example, unimodality and rank-symmetry of posets of (weak) compositions has been studied previously in \cite{sagan, zhong}. 

We now consider the poset of \emph{weak compositions} $\alpha=(a_1,\dots, a_k)$ of length $k$, denoted $\mc{C}_k$. This poset is ranked by the \emph{size} $\ab{\alpha}$ of a composition, that is $\ab{\alpha}=\sum_{i=1}^k a_i$. For $\alpha=(a_1,\dots,a_k)$ and $\beta=(b_1,\dots,b_k)$ in $\mc{C}_k$, it is also clear that $\alpha\preceq \beta$ if and only if $a_i \leq b_i$ for all $i$.  We can further restrict each part of the composition to be bounded above by some $c_i$.
\begin{lemma} \label{lem:weak_comp}
    Let $\gamma=(c_1,\dots,c_k)\in \mc{C}_k$, and denote the poset of weak compositions contained inside $\gamma$ by $\mc{C}(\gamma)$. Then:
    \begin{enumerate}
        \item $\mc{C}(\gamma)\cong \ul{c_1+1}\x \cdots \x \ul{c_k+1}$.
        \item $\mc{C}(\gamma)$ is a symmetric chain order.
    \end{enumerate}
\end{lemma}
\begin{proof}
    Item $(1)$ is immediate from the description of the order relation inherited from $\mc{C}_k$. Item $(2)$ then follows from item (1) and \Cref{lem:product_SCD}.
\end{proof} 
\begin{remark}
    As observed in \cite{sagan}, the statements of \Cref{lem:weak_comp} are rarely true for posets of ``strong" compositions, meaning those with no internal zeros, restricted to a bounding shape. As with partitions, these posets may not even be rank-symmetric, though they are rank-unimodal when restricted to a rectangle. However, strong compositions do not appear to arise in connection with Bruhat--Chevalley order as naturally as weak compositions. An explicit description of a symmetric chain decomposition of $\mc{C}(\gamma)$ for $\gamma$  a rectangle shape is given in \cite{zhong}.
\end{remark}

 Returning to the symmetric group $\mc{S}_n$ with simple reflections $S$, let $K=S-\set{s_k}$ so that $\W_{K}=\mc{S}_k \x \mc{S}_{n-k}$. Then among $\W^{K},$ we have a longest element whose one-line notation is 
\begin{equation*}
     w_1\dots w_n=(n-k+1)\ \dots \ n\ 1 \ \dots \ n-k,
\end{equation*}
and whose code is the rectangular shape $((n-k)^k)$ using the exponential notation for repeated parts of compositions. More generally, for $0\leq s \leq n-k$ we have elements of the form 
\begin{equation}
    \label{eq:w_s}
    w_{\mb{s}}=(s+1)\dots(s+k) 1 \dots s (s+k+1)\dots n    
\end{equation}
which also have rectangular codes $(s^k)$. As elements of $\W^K$, these $w_\mb{s}$ correspond to sub-Grassmannian Schubert varieties. The next result describes symmetric chain decompositions of the posets $[e,w_{\mb{s}}]^J$. 

\begin{theorem}
    \label{thm:interval_bij}
    Let $w_\mb{s}\in \W^J$ be the permutation of \eqref{eq:w_s}. Then there is an order preserving bijection given by $\code^{-1}:
    \mc{C}((s^k))\to [e, w_\mb{s}]^J$. Hence $[e,w_{\mb{s}}]^J$ is an SCO.
\end{theorem}
\begin{proof}
    By the \Cref{prop:Sn_SCD}, the map $\code^{-1}$ restricted to weak compositions inside the rectangle $(s^k)$ is order-preserving. It remains to verify that for $u\in \W^J$, $u\BC w_\mb{s}$  only if $\code(u)$ is contained in $\mc{C}((s^k))\cong(\ul{s+1})^k$; that is, the restriction of $\code^{-1}$ to $\mc{C}((s^k))$ surjects onto $[e,w_\mb{s}]^J$.

    For ease of notation, let $w_1\dots w_n$ be the one line notation of $w_\mb{s}$ for some fixed $s$, and let $u_1\dots u_n$ be the one-line notation for some $u \in \W^J$.  
    Further let $\code(u)=(a_1,\dots,a_k)$.
   
    We argue that $\code(u)\not\in (\ul{s+1})^k$ implies $u\not\BC w_{\mb{s}}$ in $\W^J$. Supposing that $\code(u)\not\in (\ul{s+1})^k$, let $i$ be the least index such that $a_i> s$. 

   Recall from \Cref{rk:code_RCM} that $a_i=u_i-r_u(u_i,i)$.  Since $i$ is the least index with $a_i>s$, 
   \[a_j=u_j-r_u(u_j,j)\leq s\quad \text{ for } 1\leq j < i.\]
    Thus,
    \begin{equation}\label{eq:ujwj}
        u_j \leq s+r_u(u_j,j)\leq s+j=w_j\quad \text{ for } 1\leq j < i.
    \end{equation}
    Next we argue that $r_u(u_i,i)=i$; that is $u_i>u_j$ for all $1\leq j < i$. Towards a contradiction, suppose otherwise that 
    \[m:=\#\set{j \mid \ 1\leq j < i,\; u_j>u_i}\geq 1.\]
    Note that then $r_u(u_i,i)=i-m$. Since by \eqref{eq:ujwj} we know  $u_j\leq s+i-1$ for $1\leq j<i$, we then have $u_i\leq s+i-1-m$. This yields
    \[a_i=u_i-r_u(u_i,i) \leq s+i-1-m-(i-m)=s-1,\]
    which contradicts that $a_i>s$. We may now assume $r_u(u_i,i)=i$.

    Now since $a_i=u_i-r_u(u_i,i)>s=w_i-r_w(w_i,i)$, we have
    \begin{equation*}
        u_i>w_i+r_u(u_i,i)-r_w(w_i,i).
    \end{equation*}
    But  $r_u(u_i,i)=r_w(w_i,i)=i$
   
    now implies $u_i>w_i$. This means that \[i=r_w(w_i,i)>r_u(w_i,i),\] and
   
    thus $r_u\not\geq r_w$ and $u\not\BC w$. This completes the proof.    
\end{proof}

\section{Remarks and future directions}\label{sec:future}
\subsection{Quotients of SCOs}
We note another approach to \Cref{prob:rect}, also mentioned in \cite[Section 4.2]{duffusThayer}. Recall that for a given poset $(\mc{P},\leq)$, a map $f:\mc{P}\to\mc{P}$ is said to be an \emph{automorphism} of $\mc{P}$ if both:
\begin{enumerate}
    \item $f$ is bijective, and
    \item $u\leq v$ if and only $f(u)\leq f(v)$.
\end{enumerate}
The automorphisms of a poset $\mc{P}$ form a group $\Aut(\mc{P})$ under composition of maps. Given a subgroup $G\sus \Aut(\mc{P})$, one can also form the quotient poset $\mc{P}/G$ defined as follows.
\begin{enumerate}
    \item The elements of $\mc{P}/G$ are the orbits of $G$, that is $[u]\in \mc{P}/G$ is the equivalence class \[[u]=\set{v\in \mc{P} \mid \exists g\in G,\ v=g\cdot u}.\]
    \item We say that $[u] \leq [v]$ if there are $u'\in [u]$ and $v'\in [v]$ such that $u'\leq v'$.
\end{enumerate}

For any $1\leq k < n$, the poset $\mc{C}((n-k)^k)=\prod_{i=1}^k \ul{n-k+1}$ admits an action by $\mc{S}_k$ which permutes the chains as automorphisms of the poset. The special case $n-k=1$ makes $\mc{C}((1)^k)$ isomorphic to the Boolean lattice $B_k$ of containment relations between subsets of a $k$-element set. 

Several authors have studied when $\mc{C}((1)^k)/G$ is a symmetric chain order; see \cite{duffusThayer} for a summary and references. That article also extends the most general known result for the Boolean poset \cite{dmt}, showing that $\mc{C}((n-k)^k)/G$ is a symmetric chain order whenever $G\sus \mc{S}_k$ is generated by disjoint cycles. An conjecture attributed to Canfield--Mason in \cite{dmt} states that $B_k/G$ admits an SCD for any subgroup of $G\sus S_k$.

\begin{prop}
    The quotient $\mc{C}((n-k)^k)/\mc{S}_k$ is isomorphic to $L(k,n-k)$.
\end{prop}
\begin{proof}
    For each orbit of $\mc{S}_k$, there is a unique representative whose parts are weakly decreasing. That is, there is a unique partition in each orbit. Two partitions $\la, \mu$ are related in $L(k,n-k)$ if and only $\la_i\geq \mu_i$ for each $1\leq i\leq k$. Since this relation is also satisfied when they are viewed as compositions, it is clear that the quotient poset is isomorphic. 
\end{proof}

Hence, as noted in \cite[Problem 2]{duffusThayer} the natural generalization of the Canfield--Mason conjecture to powers of chains would imply an affirmative answer to \Cref{prob:rect}, when taking $G=\mc{S}_k$. For more on this perspective, see also \cite{stanley1982}.

\subsection{Generalized codes and factorization of Poincar\'e polynomials} \label{sec:final.codes}
Bolognini--Sentinelli generalized ``Lehmer codes" to arbitrary Coxeter systems $(\W,S)$ as order-preserving bijections from a product of chains to $(\W, \BC)$. Observe that this notion of generalized codes actually corresponds to the inverse map $\code^{-1}$ of the usual code on $\mc{S}_n$ (\Cref{lem:staircase bijection}).

In type $B_n$, such a generalized code can be constructed as in \Cref{rk:code.BPseq} by taking sequences of sets generators $\set{J_i}$ such that the successive quotients $\W_{J_i}^{J_{i-1}}$ are chains \cite[Theorem 5.4]{bologniniSentinelli2025}. In type $D_n$, a slightly more complicated construction is required to construct a generalized code for the whole poset $(\W, \BC)$ \cite[Theorem 5.8]{bologniniSentinelli2025}. 

While not cast in the language of codes, building on work of Gasharov \cite{gasharov}, Billey had previously shown that there exists an order-preserving bijection from a product of chains to $[e,w]$ for any rationally smooth $w\in \W$ of type $A$ or $B$. From this, she also concluded that such intervals are SCOs \cite[Corollary 1.3]{billey98}. Gaetz--Gao retook the language of ``generalized Lehmer codes" from \cite{bologniniSentinelli2025} and used their results to show existence of generalized codes for $w$ rationally smooth in $\W$ of type $D_n$ \cite[Theorem  1.8]{gaetzGao}, implying the existence of SCDs for such $[e,w]$. They also consider other irreducible finite types. 

We propose the following terminology for generalizing codes to arbitrary parabolic quotients $[e,w]^J$.
\begin{defn}\label{def:gcode}
    Let $(\W,S)$ be a Coxeter system and $w\in \W^J$ for some $J\sus S$. We say that $[e,w]^J$ admits a \dfn{g-code} if there is an order- and rank-preserving bijection
    \[\operatorname{code^{-1}}: \mc{C}(\alpha)\to [e,w]^J\]    for some composition $\alpha\vDash \ell(w)$ with $\alpha\in \mc{C}_{\ab{S}}$.
\end{defn}
A g-code then automatically provides a factorization of the Poincar\'e polynomial $P_w^J(q)=\prod_{i=1}^{\ab{S}}[\alpha_i]_q$ as a product of $q$-integers. One difficulty in studying arbitrary g-codes is that they may not arise as the restriction of those on the full poset of $\W$ or $\W^J$. We mentioned an instance of this possibility in \Cref{rk:code.restrict}. Already for $w\in \W$, in the construction of a complete BP decomposition, it may be necessary to pass to $w^{-1}$ to obtain a Grassmannian BP decomposition with respect to $(J,K)$ such that $\W^K$ is a chain. Alternatively, one may need to take a left parabolic decomposition $uv=w$ to factorize a given $P_w(q)$ as the product of a $q$-integer and another palindromic Poincar\'e polynomial $P_{w'}(q)$; see \cite[Theorem 5.1]{richmondSlofstra} for a precise statement.

 We plan to study the existence of g-codes in greater detail, as well as the ``factoring'' of intervals $[e,w]^J$ into products of Grassmannian intervals as presented in \Cref{prop:RS_cor}(1), in future work.
 While a g-code provides an SCD, clearly it is not a necessary condition for the existence of an SCD. Furthermore, the factorization of the rank-generating polynomial of an arbitrary rank-symmetric poset into $q$-integers is not generally sufficient for the existence of an order-preserving bijection from a product of chains. This makes the coincidence of these properties for intervals $[e,w]$ in classical types remarkable. 

\subsection{Representation theoretic coefficients}
The question of rank-unimodality, the Peck property, and SCDs for finite posets in general and $L(k,n-k)$ in particular has historically drawn attention due to its relationship to representations of the Lie algebra $\mathfrak{sl}_2(\C)$. Stanley credits the observation of this connection (via Hard Lefschetz Theorem, as proved by Chern \cite{chern}) to Tony Iarrobino \cite[p. 174]{stanley1980}. An expository account of this perspective can be found in \cite{proctor82sl2}. More recently, a specific connection to Kronecker and plethysm coefficients has been observed, which we now describe.

Letting 
\[L_d(k,n-k):=\set{\lambda\in L(k,n-k) \mid \lambda \vdash d},\]
the construction of an SCD for $L(k,n-k)$ provides a complete matching from $L_d(k,n-k)$ to $L_{d+1}(k,n-k)$ for $0\leq d \leq \floor{k(n-k)/2}$. Thus, if $L(k,n-k)$ admits an SCD, the quantity 
\[c_d(k,n-k):=\#L_{d}(k,n-k)-\# L_{d-1}(k,n-k)\]
can be interpreted as the number of chains of the SCD which start at rank $d$. 

The question of combinatorially interpreting Kronecker and plethysm coefficients have vexed mathematicians for decades \cite{stanley2000}. The \emph{Kronecker coefficient} $g(\lambda,\mu,\nu)$ indexed by partitions $\lambda,\mu,\nu\vdash n$ gives the multiplicity of an irreducible representation $V_\nu$ of $\mc{S}_n$ in the tensor product $V_\lambda\ox V_\mu$. It was first observed by Clausen--Meier \cite{clausenMeier} that the \emph{two-row rectangular} Kronecker coefficients satisfy
\[g(((n-k)^k),((n-k)^k),(k(n-k)-d,d))=c_d(k,n-k).\]
Thus, an SCD for $L(k,n-k)$ would provide a positive combinatorial interpretation for these coefficients.

On the other hand, the plethysm coefficients $a_{\lambda \mu}^\nu$ describe the decomposition of the plethystic substitution of Schur functions,
\[s_\lambda[s_\mu]=\sum_{\nu} a_{\lambda\mu}^\nu s_\nu.\]
The expression above can be regarded as the character of a polynomial general linear group representation decomposed into irreducibles. In particular, $a_{(n-k) (k)}^\nu$ gives the expansion of $\sym^{n-k}\sym^k V$ into irreducibles where $V$ is the defining representation of the general linear group $GL(V)$. For $V\cong \C^2$, the length of $\nu$ is at most two and it is known that
\[a_{(n-k) (k)}^\nu=c_{\nu_2}(k,n-k),\]
where $\nu=(\nu_1,\nu_2)$
See, for instance, \cite[Section 5]{ossz} and references for more on direct connections between SCDs of $L(k,n-k)$ and plethysm coefficients.

A combinatorial interpretation for $c_d(k,n-k)$ was recently described in \cite{pakPanovaSwanson}, based on a recursive construction of an SCD for a modified version of $L(k,n-k)$ devised by O'Hara \cite{KOH}. Despite these advances, \Cref{prob:rect} remains open. 

\subsection*{Acknowledgements}
A. Bingham is supported by FONDECYT-ANID grant 3250472.

\subsection*{Statement on AI use} We used a mix of large-language models for some light proofreading, grammatical help, and assistance making some of the diagrams in~\Cref{sec:results}.  We did not use any AI for exploration, brainstorming, computation, mathematical reasoning, proving, or in any other capacity than those mentioned.

\bibliography{bio-scd-WP} 

@book{anderson,
 author = {Anderson, Ian},
 title = {Combinatorics of finite sets.},
 edition = {Corrected reprint of the 1989 edition},
 isbn = {0-486-42257-7},
 year = {2002},
 publisher = {Mineola, NY: Dover Publications},
 language = {English},
 zbMATH = {1792287},
 Zbl = {0995.05001}
}

@article{billey98,
 author = {Billey, Sara C.},
 title = {Pattern avoidance and rational smoothness of {Schubert} varieties},
 fjournal = {Advances in Mathematics},
 journal = {Adv. Math.},
 issn = {0001-8708},
 volume = {139},
 number = {1},
 pages = {141--156},
 year = {1998},
 language = {English},
 doi = {10.1006/aima.1998.1744},
 zbMATH = {1236157},
 Zbl = {0936.14036}
}

@book{billeyLakshmibai,
 author = {Billey, Sara and Lakshmibai, V.},
 title = {Singular loci of {Schubert} varieties},
 fseries = {Progress in Mathematics},
 series = {Prog. Math.},
 issn = {0743-1643},
 volume = {182},
 isbn = {3-7643-4092-4},
 year = {2000},
 publisher = {Boston, MA: Birkh{\"a}user},
 language = {English},
 zbMATH = {1520629},
 Zbl = {0959.14032}
}

@article{billeyPostnikov,
 author = {Billey, Sara and Postnikov, Alexander},
 title = {Smoothness of {Schubert} varieties via patterns in root subsystems},
 fjournal = {Advances in Applied Mathematics},
 journal = {Adv. Appl. Math.},
 issn = {0196-8858},
 volume = {34},
 number = {3},
 pages = {447--466},
 year = {2005},
 language = {English},
 doi = {10.1016/j.aam.2004.08.003},
 zbMATH = {2212819},
 Zbl = {1072.14065}
}

@book{bjornerBrenti,
 author = {Bj{\"o}rner, Anders and Brenti, Francesco},
 title = {Combinatorics of {Coxeter} groups},
 fseries = {Graduate Texts in Mathematics},
 series = {Grad. Texts Math.},
 issn = {0072-5285},
 volume = {231},
 isbn = {3-540-44238-3},
 year = {2005},
 publisher = {New York, NY: Springer},
 language = {English},
 doi = {10.1007/3-540-27596-7},
 zbMATH = {2186170},
 Zbl = {1110.05001}
}

@misc{bologniniSentinelli2025,
  author        = {Bolognini, Davide and Sentinelli, Paolo},
  title         = {The {Lehmer} complex of a {Bruhat} interval},
  year          = {2026},
  howpublished  = {Preprint, arXiv:2501.03037 [math.CO], \url{https://arxiv.org/abs/2501.03037}},
  eprint        = {2501.03037},
  archivePrefix = {arXiv},
  primaryClass  = {math.CO},
}

@book{borelLAG,
 author = {Borel, Armand},
 title = {Linear algebraic groups.},
 edition = {2nd enlarged ed.},
 fseries = {Graduate Texts in Mathematics},
 series = {Grad. Texts Math.},
 issn = {0072-5285},
 volume = {126},
 isbn = {0-387-97370-2},
 year = {1991},
 publisher = {New York etc.: Springer-Verlag},
 language = {English},
 zbMATH = {50185},
 Zbl = {0726.20030}
}

@article{carrellKuttler,
 author = {Carrell, James B. and Kuttler, Jochen},
 title = {Smooth points of {{\(T\)}}-stable varieties in {{\(G\)}}/{{\(B\)}} and the {Peterson} map},
 fjournal = {Inventiones Mathematicae},
 journal = {Invent. Math.},
 issn = {0020-9910},
 volume = {151},
 number = {2},
 pages = {353--379},
 year = {2003},
 language = {English},
 doi = {10.1007/s00222-002-0256-5},
 zbMATH = {1965467},
 Zbl = {1041.14020}
}

@incollection{carrellPeterson,
 author = {Carrell, James B.},
 title = {The {Bruhat} graph of a {Coxeter} group, a conjecture of {Deodhar}, and rational smoothness of {Schubert} varieties},
 booktitle = {Algebraic groups and their generalizations: classical methods. Summer Research Institute on algebraic groups and their generalizations, July 6- 26, 1991, Pennsylvania State University, University Park, PA, USA},
 isbn = {0-8218-1540-7; 0-8218-1497-4},
 pages = {53--61},
 year = {1994},
 publisher = {Providence, RI: American Mathematical Society},
 language = {English},
 zbMATH = {647429},
 Zbl = {0818.14020}
}

@article{chern,
 author = {Chern, Shiing-shen},
 title = {On a generalization of {Kaehler} geometry},
 fjournal = {Princeton Mathematical Series},
 journal = {Princeton Math. Ser.},
 volume = {12},
 pages = {103--121},
 year = {1957},
 language = {English},
 zbMATH = {3127980},
 Zbl = {0078.14103}
}

@article{debruijn_etal,
 author = {de Bruijn, N. G. and van Ebbenhorst Tengbergen, C. and Kruyswijk, D.},
 title = {On the set of divisors of a number},
 fjournal = {Nieuw Archief voor Wiskunde. Tweede Serie},
 journal = {Nieuw Arch. Wiskd., II. Ser.},
 issn = {0028-9825},
 volume = {23},
 pages = {191--193},
 year = {1951},
 language = {English},
 zbMATH = {3065933},
 Zbl = {0043.04301}
}

@article{clausenMeier,
 author = {Clausen, Michael and Meier, Helga},
 title = {Extreme irreducible constituents in tensor representations of symmetric groups},
 fjournal = {Bayreuther Mathematische Schriften},
 journal = {Bayreuther Math. Schr.},
 issn = {0172-1062},
 volume = {45},
 pages = {1--17},
 year = {1993},
 language = {English},
 zbMATH = {553948},
 Zbl = {0799.05004}
}

@article{deodhar85,
 author = {Deodhar, Vinay V.},
 title = {Local {Poincar{\'e}} duality and non-singularity of {Schubert} varieties},
 fjournal = {Communications in Algebra},
 journal = {Commun. Algebra},
 issn = {0092-7872},
 volume = {13},
 pages = {1379--1388},
 year = {1985},
 language = {English},
 doi = {10.1080/00927878508823227},
 zbMATH = {3927108},
 Zbl = {0579.14046}
}

@article{dmt,
 author = {Duffus, Dwight and McKibben-Sanders, Jeremy and Thayer, Kyle},
 title = {Some quotients of chain products are symmetric chain orders},
 fjournal = {The Electronic Journal of Combinatorics},
 journal = {Electron. J. Comb.},
 issn = {1077-8926},
 volume = {19},
 number = {2},
 pages = {11},
 note = {Id/No p46},
 year = {2012},
 language = {English},
 url = {www.combinatorics.org/ojs/index.php/eljc/article/view/v19i2p46},
 zbMATH = {6098333},
 Zbl = {1253.06002}
}

@article{duffusThayer,
 author = {Duffus, Dwight and Thayer, Kyle},
 title = {Symmetric chain decompositions of quotients by wreath products.},
 fjournal = {The Electronic Journal of Combinatorics},
 journal = {Electron. J. Comb.},
 issn = {1077-8926},
 volume = {22},
 number = {2},
 pages = {15},
 note = {Id/No p2.35},
 year = {2015},
 language = {English},
 url = {www.combinatorics.org/ojs/index.php/eljc/article/view/v22i2p35},
 zbMATH = {6441864},
 Zbl = {1315.06003}
}

@book{engel,
 author = {Engel, Konrad},
 title = {Sperner theory},
 fseries = {Encyclopedia of Mathematics and Its Applications},
 series = {Encycl. Math. Appl.},
 issn = {0953-4806},
 volume = {65},
 isbn = {0-521-45206-6},
 year = {1997},
 publisher = {Cambridge: Cambridge University Press},
 language = {English},
 zbMATH = {1016362},
 Zbl = {0868.05001}
}

@misc{gaetzGao,
  author        = {Gaetz, Christian and Gao, Yibo},
  title         = {{Billey--Postnikov} posets, rationally smooth {Schubert} varieties, and {Poincar{\'e}} duality},
  year          = {2026},
  howpublished  = {Preprint, arXiv:2512.08168 [math.CO], \url{https://arxiv.org/abs/2512.08168}},
  eprint        = {2512.08168},
  archivePrefix = {arXiv},
  primaryClass  = {math.CO},
}

@article{gasharov,
 author = {Gasharov, Vesselin},
 title = {Factoring the {Poincar{\'e}} polynomials for the {Bruhat} order on {{\(S_n\)}}},
 fjournal = {Journal of Combinatorial Theory. Series A},
 journal = {J. Comb. Theory, Ser. A},
 issn = {0097-3165},
 volume = {83},
 number = {1},
 pages = {159--164},
 year = {1998},
 language = {English},
 doi = {10.1006/jcta.1997.2861},
 zbMATH = {1212920},
 Zbl = {0935.05087}
}

@book{Humphreys,
 author = {Humphreys, James E.},
 title = {Reflection groups and {Coxeter} groups},
 fseries = {Cambridge Studies in Advanced Mathematics},
 series = {Camb. Stud. Adv. Math.},
 volume = {29},
 isbn = {0-521-37510-X},
 year = {1990},
 publisher = {Cambridge etc.: Cambridge University Press},
 language = {English},
 zbMATH = {47598},
 Zbl = {0725.20028}
}

@book{humphreysLAG,
 author = {Humphreys, James E.},
 title = {Linear algebraic groups. {Corr}. 2nd printing},
 fseries = {Graduate Texts in Mathematics},
 series = {Grad. Texts Math.},
 issn = {0072-5285},
 volume = {21},
 year = {1981},
 publisher = {Springer, Cham},
 language = {English},
 zbMATH = {3737912},
 Zbl = {0471.20029}
}

@book{kumar,
 author = {Kumar, Shrawan},
 title = {Kac-{Moody} groups, their flag varieties and representation theory},
 fseries = {Progress in Mathematics},
 series = {Prog. Math.},
 issn = {0743-1643},
 volume = {204},
 isbn = {0-8176-4227-7},
 year = {2002},
 publisher = {Boston, MA: Birkh{\"a}user},
 language = {English},
 zbMATH = {1821146},
 Zbl = {1026.17030}
}

@article{lambertRabelo2026,
 author = {Lambert, Jordan and Rabelo, Lonardo},
 title = {The {Bruhat} order on symmetric groups via intrinsic coverings of compositions},
 fjournal = {The Australasian Journal of Combinatorics},
 journal = {Australas. J. Comb.},
 issn = {1034-4942},
 volume = {95},
 pages = {257--282},
 year = {2026},
 language = {English},
 url = {ajc.maths.uq.edu.au/pdf/95/ajc_v95_p257.pdf},
 zbMATH = {8234875}
}

@article{laisant,
 author = {Laisant, C. A.},
 title = {Sur la num{\'e}ration factorielle, application aux permutations.},
 fjournal = {Bulletin de la Soci{\'e}t{\'e} Math{\'e}matique de France},
 journal = {Bull. Soc. Math. Fr.},
 issn = {0037-9484},
 volume = {16},
 pages = {176--183},
 year = {1888},
 language = {French},
 doi = {10.24033/bsmf.378},
 url = {https://eudml.org/doc/85654},
 zbMATH = {2693718},
 JFM = {20.0175.03}
}

@article{lascoux,
 author = {Lascoux, Alain},
 title = {Ordering the symmetric group: why using the {Iwahori}-{Hecke} algebra?},
 fjournal = {Documenta Mathematica},
 journal = {Doc. Math.},
 issn = {1431-0635},
 volume = {Extra Vol.},
 pages = {355--364},
 year = {1998},
 language = {French},
 doi = {10.4171/dms/1-3/35},
 url = {https://eudml.org/doc/228538},
 zbMATH = {1184374},
 Zbl = {0909.05048}
}

@misc{lehmer,
 author = {Lehmer, D. H.},
 title = {Teaching combinatorial tricks to a computer},
 year = {1960},
 language = {English},
 howpublished = {Proc. {Symp}. {Appl}. {Math}. 10, 179-193.},
 zbMATH = {3156373},
 Zbl = {0096.00504}
}

@book{manivel,
 author = {Manivel, Laurent},
 title = {Symmetric functions, {Schubert} polynomials and degeneracy loci. {Transl}. from the {French} by {John} {R}. {Swallow}},
 fseries = {SMF/AMS Texts and Monographs},
 series = {SMF/AMS Texts Monogr.},
 issn = {1525-2302},
 volume = {6},
 isbn = {0-8218-2154-7},
 year = {2001},
 publisher = {Providence, RI: American Mathematical Society (AMS)},
 language = {English},
 zbMATH = {1678993},
 Zbl = {0998.14023}
}

@misc{ohRichmond,
 author = {Suho Oh and Edward Richmond},
 title = {Coxeter groups and {Billey}-{Postnikov} decompositions},
 year = {2024},
 howpublished = {Preprint, {arXiv}:2409.03096 [math.{CO}]},
 url = {https://arxiv.org/abs/2409.03096},
 arXiv = {arXiv:2409.03096}
}

@incollection{ohYoo,
 author = {Oh, Suho and Yoo, Hwanchul},
 title = {Bruhat order, rationally smooth {Schubert} varieties, and hyperplane arrangements},
 booktitle = {Proceedings of the 22nd annual international conference on formal power series and algebraic combinatorics, FPSAC 2010, San Francisco, USA, August 2--6, 2010},
 pages = {965--972},
 year = {2010},
 publisher = {Nancy: The Association. Discrete Mathematics \& Theoretical Computer Science (DMTCS)},
 language = {English},
 url = {www.dmtcs.org/dmtcs-ojs/index.php/proceedings/article/view/dmAN0173},
 zbMATH = {6806847},
 Zbl = {1374.05237}
}

@article{OSSZ,
 author = {Orellana, Rosa and Saliola, Franco and Schilling, Anne and Zabrocki, Mike},
 title = {From quasi-symmetric to {Schur} expansions with applications to symmetric chain decompositions and plethysm},
 fjournal = {The Electronic Journal of Combinatorics},
 journal = {Electron. J. Comb.},
 issn = {1077-8926},
 volume = {31},
 number = {4},
 pages = {29},
 note = {Id/No p4.23},
 year = {2024},
 language = {English},
 doi = {10.37236/12950},
 zbMATH = {7940966},
 Zbl = {1551.05412}
}

@article{KOH,
 author = {O'Hara, Kathleen M.},
 title = {Unimodality of {Gaussian} coefficients: {A} constructive proof},
 fjournal = {Journal of Combinatorial Theory. Series A},
 journal = {J. Comb. Theory, Ser. A},
 issn = {0097-3165},
 volume = {53},
 number = {1},
 pages = {29--52},
 year = {1990},
 language = {English},
 doi = {10.1016/0097-3165(90)90018-R},
 zbMATH = {4142040},
 Zbl = {0697.05002}
}

@misc{pakPanovaSwanson,
 author = {Igor Pak and Greta Panova and Joshua P. Swanson},
 title = {A combinatorial interpretation for certain plethysm and {Kronecker} coefficients},
 year = {2026},
 howpublished = {Preprint, {arXiv}:2511.02312 [math.{CO}]},
 url = {https://arxiv.org/abs/2511.02312},
 arXiv = {arXiv:2511.02312}
}

@article{peck,
 author = {Peck, G. W.},
 title = {Kleitman and combinatorics: {A} celebration},
 fjournal = {Discrete Mathematics},
 journal = {Discrete Math.},
 issn = {0012-365X},
 volume = {257},
 number = {2-3},
 pages = {193--224},
 year = {2002},
 language = {English},
 doi = {10.1016/S0012-365X(02)00595-2},
 zbMATH = {1838988},
 Zbl = {1027.01016}
}

@article{proctor82,
 author = {Proctor, Robert A.},
 title = {Classical {Bruhat} orders and lexicographic shellability},
 fjournal = {Journal of Algebra},
 journal = {J. Algebra},
 issn = {0021-8693},
 volume = {77},
 pages = {104--126},
 year = {1982},
 language = {English},
 doi = {10.1016/0021-8693(82)90280-0},
 zbMATH = {3763992},
 Zbl = {0486.06002}
}

@article{proctor82sl2,
 author = {Proctor, Robert A.},
 title = {Representations of sl(2,{C}) on posets and the {Sperner} property},
 fjournal = {SIAM Journal on Algebraic and Discrete Methods},
 journal = {SIAM J. Algebraic Discrete Methods},
 issn = {0196-5212},
 volume = {3},
 pages = {275--280},
 year = {1982},
 language = {English},
 doi = {10.1137/0603026},
 zbMATH = {3781238},
 Zbl = {0496.06004}
}

@article{proctorSaksSturtevant,
 author = {Proctor, Robert A. and Saks, Michael E. and Sturtevant, Dean G.},
 title = {Product partial orders with the {Sperner} property},
 fjournal = {Discrete Mathematics},
 journal = {Discrete Math.},
 issn = {0012-365X},
 volume = {30},
 pages = {173--180},
 year = {1980},
 language = {English},
 doi = {10.1016/0012-365X(80)90118-1},
 zbMATH = {3717374},
 Zbl = {0458.06001}
}

@article{richmondSlofstra,
 author = {Richmond, Edward and Slofstra, William},
 title = {Billey-{Postnikov} decompositions and the fibre bundle structure of {Schubert} varieties},
 fjournal = {Mathematische Annalen},
 journal = {Math. Ann.},
 issn = {0025-5831},
 volume = {366},
 number = {1-2},
 pages = {31--55},
 year = {2016},
 language = {English},
 doi = {10.1007/s00208-015-1299-4},
 zbMATH = {6643569},
 Zbl = {1383.14013}
}

@article{riess,
  author  = {Riess, W.},
  title   = {Zwei {Optimierungsprobleme} auf {Ordnungen}},
  journal = {Arbeitsberichte des Instituts f{\"u}r mathematische Maschinen und Datenverarbeitung (Informatik), Erlangen},
  volume  = {Bank II},
  number  = {5},
  pages   = {50--57},
  year    = {1978},
}

@article{ryan,
 author = {Ryan, Kevin M.},
 title = {On {Schubert} varieties in the flag manifold of {Sl}(n,{{\({\mathbb C})\)}}},
 fjournal = {Mathematische Annalen},
 journal = {Math. Ann.},
 issn = {0025-5831},
 volume = {276},
 pages = {205--224},
 year = {1987},
 language = {English},
 doi = {10.1007/BF01450738},
 url = {https://eudml.org/doc/164193},
 zbMATH = {3927107},
 Zbl = {0579.14045}
}

@article{sagan,
 author = {Sagan, Bruce E.},
 title = {Compositions inside a rectangle and unimodality},
 fjournal = {Journal of Algebraic Combinatorics},
 journal = {J. Algebr. Comb.},
 issn = {0925-9899},
 volume = {29},
 number = {4},
 pages = {405--411},
 year = {2009},
 language = {English},
 doi = {10.1007/s10801-008-0140-5},
 zbMATH = {5607590},
 Zbl = {1226.05262}
}

@article{sperner,
 author = {Sperner, E.},
 title = {Ein {Satz} {\"u}ber {Untermengen} einer endlichen {Menge}.},
 fjournal = {Mathematische Zeitschrift},
 journal = {Math. Z.},
 issn = {0025-5874},
 volume = {27},
 pages = {544--548},
 year = {1928},
 language = {German},
 doi = {10.1007/BF01171114},
 url = {https://eudml.org/doc/167993},
 zbMATH = {2575694},
 JFM = {54.0090.06}
}

@book{springerLAG,
 author = {Springer, T. A.},
 title = {Linear algebraic groups.},
 edition = {Reprint of the 1998 2nd ed.},
 fseries = {Modern Birkh{\"a}user Classics},
 series = {Mod. Birkh{\"a}user Class.},
 issn = {2197-1803},
 isbn = {978-0-8176-4839-8},
 year = {2009},
 publisher = {Basel: Birkh{\"a}user},
 language = {English},
 zbMATH = {5359373},
 Zbl = {1202.20048}
}

@article{stanley1980,
 author = {Stanley, Richard P.},
 title = {Weyl groups, the hard {Lefschetz} theorem, and the {Sperner} property},
 fjournal = {SIAM Journal on Algebraic and Discrete Methods},
 journal = {SIAM J. Algebraic Discrete Methods},
 issn = {0196-5212},
 volume = {1},
 pages = {168--184},
 year = {1980},
 language = {English},
 doi = {10.1137/0601021},
 zbMATH = {3790431},
 Zbl = {0502.05004}
}

@article{stanley1982,
 author = {Stanley, Richard P.},
 title = {Some aspects of groups acting on finite posets},
 fjournal = {Journal of Combinatorial Theory. Series A},
 journal = {J. Comb. Theory, Ser. A},
 issn = {0097-3165},
 volume = {32},
 pages = {132--161},
 year = {1982},
 language = {English},
 doi = {10.1016/0097-3165(82)90017-6},
 zbMATH = {3781235},
 Zbl = {0496.06001}
}

@incollection{stanley2000,
 author = {Stanley, Richard P.},
 title = {Positivity problems and conjectures in algebraic combinatorics},
 booktitle = {Mathematics: frontiers and perspectives},
 isbn = {0-8218-2070-2},
 pages = {295--319},
 year = {2000},
 publisher = {Providence, RI: American Mathematical Society (AMS)},
 language = {English},
 zbMATH = {1503622},
 Zbl = {0955.05111}
}

@book{stanleyEC1,
 author = {Stanley, Richard P.},
 title = {Enumerative combinatorics. {Vol}. 1.},
 edition = {2nd ed.},
 fseries = {Cambridge Studies in Advanced Mathematics},
 series = {Camb. Stud. Adv. Math.},
 volume = {49},
 isbn = {978-1-107-60262-5; 978-1-107-01542-5; 978-1-139-20056-1},
 year = {2012},
 publisher = {Cambridge: Cambridge University Press},
 language = {English},
 url = {www.cambridge.org/de/knowledge/isbn/item6832283/?site_locale=de_DE},
 zbMATH = {6016068},
 Zbl = {1247.05003}
}

@article{sylvester,
 author = {Sylvester, J. J.},
 title = {Proof of the hitherto undemonstrated fundamental theorem of invariants.},
 fjournal = {The London, Edinburgh, and Dublin Philosophical Magazine and Journal of Science, V. Series},
 journal = {Phil. Mag. (5)},
 issn = {0031-8086},
 volume = {5},
 pages = {178--188},
 year = {1878},
 language = {English},
 zbMATH = {2710739},
 JFM = {10.0082.02}
}

@article{wen24,
 author = {Wen, Xiangdong},
 title = {A symmetric chain decomposition of {{\(L(5,n)\)}}},
 fjournal = {Enumerative Combinatorics and Applications},
 journal = {Enumer. Comb. Appl.},
 issn = {2710-2335},
 volume = {4},
 number = {1},
 pages = {8},
 note = {Id/No s2r5},
 year = {2024},
 language = {English},
 doi = {10.54550/ECA2024V4S1R5},
 zbMATH = {7964427}
}

@misc{wen26,
  author        = {Wen, Xiangdong},
  title         = {A symmetric chain decomposition of {L}(6,n)},
  year          = {2026},
  howpublished  = {Preprint, arXiv:2604.23042 [math.CO], \url{https://arxiv.org/abs/2604.23042}},
  eprint        = {2604.23042},
  archivePrefix = {arXiv},
  primaryClass  = {math.CO},
}

@article{wolper,
 author = {Wolper, James S.},
 title = {A combinatorial approach to the singularities of {Schubert} varieties},
 fjournal = {Advances in Mathematics},
 journal = {Adv. Math.},
 issn = {0001-8708},
 volume = {76},
 number = {2},
 pages = {184--193},
 year = {1989},
 language = {English},
 doi = {10.1016/0001-8708(89)90048-0},
 zbMATH = {4156654},
 Zbl = {0705.14048}
}

@article{zhong,
 author = {Zhong, Yueming},
 title = {A combinatorial proof of the unimodality and symmetry of weak composition rank sequences},
 fjournal = {Annals of Combinatorics},
 journal = {Ann. Comb.},
 issn = {0218-0006},
 volume = {27},
 number = {2},
 pages = {281--295},
 year = {2023},
 language = {English},
 doi = {10.1007/s00026-022-00624-0},
 zbMATH = {7711276},
 Zbl = {1518.05012}
}
\bibliographystyle{alpha}

\end{document}